\documentclass[12pt]{article}
\usepackage{amsmath,amssymb,latexsym, epsfig}
\title{Tunnel effect for Kramers-Fokker-Planck type operators: return to equilibrium and applications}
\author{Fr\'ed\'eric H\'erau\\\small Laboratoire de
Math\'ematiques\\\small Universit\'e de Reims\\\small Moulin de la Housse B.P.
1039\\\small 51687 Reims cedex 2, France\\\small and UMR 6056 CNRS\\\small herau@univ-reims.fr \and Michael Hitrik
\\\small Department of Mathematics
\\\small University of California \\\small Los Angeles
\\\small CA 90095-1555, USA\\\small hitrik@math.ucla.edu \and
Johannes Sj\"ostrand\\\small CMLS\\\small Ecole Polytechnique\\\small 91128 Palaiseau cedex, France \\\small and
UMR 7640 CNRS
\\\small johannes@math.polytechnique.fr}

\date{}
\def\trans{^t\hskip -2pt}

\catcode`\é = \active\def é{\'e} \catcode`\è = \active\def è{\`e}
\catcode`\ç = \active\def ç{\c{c}} \catcode`\à = \active\def
à{\`a} \catcode`\ê = \active \def ê{\^e} \catcode`\ù = \active
\def ù{\`u} \catcode`\ô = \active \def ô{\^o} \catcode`\û =
\active \def û{\^u} \catcode`\î = \active \def î{\^{\i}}
\catcode`\â = \active \def â{\^a} \catcode`\ö = \active \def
ö{{\"o}}

\newcommand{\real}{\mbox{\bf R}}

\def\abs#1{\left|#1\right|}
\def\begeq{\begin{equation}}
\def\endeq{\end{equation}}

\def\ccc{{\cal C}}
\def\fff{{\cal F}}
 \def\hhh{{\cal H}}

\def\mmm{{\cal M}} \def\ooo{{\cal O}}

\def\norm#1{\left\Vert#1\right\Vert}

\def\abs#1{\left\vert#1\right\vert}
\def\set#1{\left\{#1\right\}}

\def\sep#1{\left(#1\right)}

\newcommand{\preuve}[1][\!\!]{\bigskip\noindent{\bf Proof #1. \ \ }}
\def\fin{\hfill$\Box$\\}
\def\Ref#1{(\ref{#1})}
\def\wrtext#1{\relax\ifmmode{\leavevmode\hbox{#1}}\else{#1}\fi}

\newcommand{\eps}{\varepsilon}
\def\part#1{\frac{\partial}{\partial #1}}

\def\norm#1{||\,#1\,||}

\def\D{\partial}
\def\eps{\varepsilon}
\def\norm#1{\left\Vert#1\right\Vert}

\def\abs#1{\left\vert#1\right\vert}
\def\set#1{\left\{#1\right\}}

\def\sep#1{\left(#1\right)}

\def\Re{{\mathrm Re\,}} 

\def\defegal{\stackrel{\text{\rm def}}{=}}

\def\trans{^t\hskip -2pt}
\def\sm{\setminus }

\newcommand{\comp}{\mbox{\bf C}}

\newcommand{\nat}{\mbox{\bf N}}

\renewcommand{\Re}{\mbox{\rm Re\,}}
\renewcommand{\exp}{\mbox{\rm exp\,}}

\newcommand{\neigh}{neighborhood}

\def\Ref#1{(\ref{#1})}

\newtheorem{dref}{Definition}[section]
\newtheorem{lemma}[dref]{Lemma}
\newtheorem{theo}[dref]{Theorem}
\newtheorem{prop}[dref]{Proposition}

\renewcommand{\thedref}{\thesection.\arabic{dref}}
\newenvironment{proof}{\vspace{.3cm}\noindent{{\em Proof:}}}{\hfill$\Box$
\vspace{.2cm}}
\def\remark{\refstepcounter{dref}\bigskip\noindent\bf Remark \thedref.\rm\ }

\def\D{\partial}
\def\R{\mathbb R}\def\C{\mathbb C}\def\N{\mathbb N}
\def\E{\mathbb E}
\def\nat{\N} \def\comp{\C} \def\real{\R}

\begin{document}
\maketitle

\vspace*{1cm}
\noindent
{\bf Abstract}: In the first part of this work, we consider second order supersymmetric differential
operators in the semiclassical limit, including the Kramers-Fokker-Planck operator, such that the
exponent of the associated Maxwellian $\phi$ is a Morse function with
two local minima and one saddle point. Under suitable
additional assumptions of dynamical nature, we establish the
long time convergence to the equilibrium for the associated heat semigroup, with
the rate given by the first non-vanishing, exponentially small,
eigenvalue. In the second part of the paper, we consider the case when
the function $\phi$ has precisely one local minimum and one saddle point. We
also discuss further examples of supersymmetric operators, including
the Witten Laplacian and the infinitesimal generator for the
time evolution of a chain of classical anharmonic oscillators.


\tableofcontents
\section{Introduction and statement of the main result}
\setcounter{equation}{0}
The principal purpose of the present paper is to apply the spectral results
of~\cite{HeHiSj} to obtain a precise information concerning the large
time behavior of the heat semigroup generated by the semiclassical Kramers-Fokker-Planck operator
\begeq
\label{eq01}
P=y\cdot h\partial_x-V'(x)\cdot
h\partial_y+\frac{\gamma}{2}\left(-h\partial_y+y\right)\cdot
\left(h\partial_y+y\right),\quad x,y\in \real^n,\quad \gamma>0.
\endeq
In fact, as in~\cite{HeHiSj}, our main result will be valid for a large class of supersymmetric second order
differential operators, including (\ref{eq01}). Physically, the semiclassical limit $h\rightarrow 0$ in (\ref{eq01})
corresponds to the r\'egime of low temperatures. Recall that by supersymmetry, we mean the fact that the Kramers-Fokker-Planck
operator $P$ can be viewed as a Witten Laplacian in degree 0 associated to
a certain non-semidefinite scalar product in the spaces of
differential forms. These relations with  the Witten complex \cite{Wi82} were exhibited in the
works of J.~Tailleur, S.~Tanase-Nicola, J.~Kurchan \cite{TaTaKu} and
J.~M.~Bismut \cite{Bi},  using respectively the languages of
supersymmetry and differential forms (See also \cite{Le}
for a quick introduction to the differential form version of Bismut and \cite{BiLe}).

\medskip
The paper~\cite{HeHiSj}, which is a natural continuation of~\cite{HeSjSt}, analyzed resolvent estimates and
the low lying eigenvalues of $P$, assuming that the potential $V$ in (\ref{eq01}) is a smooth real valued
Morse function on $\real^n$, such that
\begeq
\label{eq011}
\partial^{\alpha}V ={\cal O}(1),\quad \abs{\alpha}\geq 2,
\endeq
and with
\begeq
\label{eq012}
\abs{\nabla V}\geq 1/C,\quad \wrtext{for}\quad \abs{x}\geq C>0.
\endeq
Assuming furthermore that $V$ has precisely three critical points: two local
minima, $x_{\pm 1}$, and one critical point $x_0$ of index $1$, it was established
in~\cite{HeHiSj} that for $C>0$ large enough and $h>0$ sufficiently small, the operator $P$ has precisely two
eigenvalues in the disc $D(0,h/C)=\{z\in \comp; \abs{z}<\frac{h}{C}\}$, $\mu_0$ and $\mu_1$, such that
$\mu_0=0$ and $\mu_1$ is real and of the form
\begeq
\label{eq02}
\mu_1=h\left(a_1(h) e^{-2S_1/h}+a_{-1}(h) e^{-2S_{-1}/h}\right), \quad
S_j=V(x_0)-V(x_j).
\endeq
Here $a_j$ are real with
\begeq
\label{eq03}
a_j(h)\sim a_{j,0}+ha_{j,1}+\ldots\,\quad a_{j,0}>0.
\endeq
Notice that the eigenfunction corresponding to the eigenvalue $\mu_0=0$ is the Max\-wel\-lian
\begeq
\label{eq031}
\exp\left(-\phi/h\right)\in L^2(\real^{2n}_{x,y}),\quad \phi(x,y)=\frac{y^2}{2}+V(x).
\endeq

\medskip
In the case when $V\rightarrow +\infty$ as $x\rightarrow \infty$ and $V$ has precisely one local minimum,
it follows from the results of~\cite{HeSjSt} that in a disc $D(0,Ch)$, $C\geq 1$, apart from the eigenvalue $\mu_0=0$,
the real part of the other eigenvalues is $\geq h/C$. In this case, precise results describing
the behavior of the semigroup $\exp(-tP/h)$ for large $t$, were obtained in~\cite{HeSjSt} --- the rate of the
return to equilibrium, given by the projection onto the ground state (\ref{eq031}), is then of the order
of magnitude $1$. In this work, we shall complement this study by analyzing the
question of a return to equilibrium in the presence of exponentially small eigenvalues, due to the tunneling between
the local minima.

\medskip
The following is the main result of this work, specialized to the case of the Kramers-Fokker-Planck operator
(\ref{eq01}). Here we shall also write $P$ for the m-accretive realization of the operator (\ref{eq01}) --- see also
section 2 and section 3 in~\cite{HeHiSj}.

\begin{theo}
Assume that $V$ in {\rm (\ref{eq01})} is a $C^{\infty}$ real valued Morse function satisfying {\rm (\ref{eq011})} and
{\rm (\ref{eq012})} and having precisely {\rm 3} critical points: {\rm 2} local minima, $x_{\pm 1}$, and
a critical point of index {\rm 1}, so that the disc $D(0,h/C)$ for $C>0$ large enough, contains precisely {\rm 2}
eigenvalues of $P$, $\mu_0=0$ and $\mu_1$ given in {\rm (\ref{eq02})}. Let $\Pi_j$ be the spectral projection
associated with the eigenvalue $\mu_j$, $j=0,1$. Then we have
\begeq
\label{eq04}
\Pi_j={\cal O}(1),\quad h\rightarrow 0.
\endeq
We have furthermore, uniformly as $t\geq 0$ and $h\rightarrow 0$,
\begeq
\label{eq05}
e^{-tP/h}=\Pi_0+e^{-t\mu_1/h}\Pi_1+{\cal O}(1) e^{-t/C}, \quad C>0,\quad \wrtext{in}\quad {\cal L}(L^2,L^2).
\endeq
\end{theo}

\medskip
\noindent
The structure of the article is as follows: In section 2, relying upon the results of~\cite{HeHiSj},~\cite{HeSjSt},
we establish the basic a priori coercivity estimate for the operator $P$ in a suitable exponentially weighted space, introduced in
~\cite{HeHiSj}. In this sense it can be interpreted as an hypocoercive estimate (see e.g. \cite{Her07}, \cite{Vil07}). In section 3 it is then quite straightforward to prove Theorem 1.1 in its general form,
by combining the results of section 2 together with the analysis of~\cite{HeHiSj}.

The second part of the paper, consisting of sections 4--6 is of a somewhat different nature, complementing and extending
the previous analysis. In section 4, we study the case
left out in~\cite{HeHiSj}, when the potential $V$ in (\ref{eq01}) has precisely two
critical points: one local minimum and a critical point of index one. In this case, $0$ is not an eigenvalue of $P$, and
the large time behavior of the heat semigroup is governed entirely by the first, exponentially small, eigenvalue. In
section 5, we give some examples and describe a probabilistic framework in which the Witten Laplacian and the
Kramers-Fokker-Planck operator both arise naturally. Finally, in section 6, we describe another example of a
supersymmetric operator, for which the question of a convergence to equilibrium is of interest, namely a chain of
classical interacting anharmonic oscillators, coupled to a heat bath. We show how to adapt the analysis
of~\cite{HeHiSj} to cover also this case.

\medskip
\noindent
{\bf Acknowledgments:} A part of this work was carried out in May of 2007, when the second author was visiting
Universit\'e Paris 13 being on leave from UCLA. It is a great pleasure
for him to thank its D\'epartement de Math\'ematiques, and in particular
Alain Grigis, for the extraordinary hospitality and excellent working
conditions. His best thanks are also due to the Centre de
Math\'ematiques of \'Ecole Polytechnique for the generous hospitality
in June of 2007. The partial support
of his research by the National Science Foundation under grant DMS-0653275 and by an Alfred P. Sloan Research Fellowship
is also gratefully acknowledged.

\section{An (hypo-)coercive estimate}
\setcounter{equation}{0}
The purpose of this section is to establish an a priori estimate for $P$, instrumental in proving Theorem 1.1.
This estimate will imply the exponential decay for the heat semigroup $\exp(-tP/h)$ in $L^2$,
when restricted to the kernel of the spectral projection corresponding to the eigenvalues $\mu_0=0$ and $\mu_1$ in
(\ref{eq02}). When doing so, as in~\cite{HeHiSj}, rather than working directly with (\ref{eq01}), we shall consider a
broader class of scalar real second
order non-elliptic non-selfadjoint operators on $\real^n$. For completeness, we shall now recall, following~\cite{HeHiSj},
the definition and the main assumptions concerning this class.

\medskip
Let us consider
\begin{eqnarray}
\label{eq1.1}
\hskip -0.4cm
P & = & \sum_{j,k=1}^n hD_{x_j}\circ b_{j,k}(x)\circ hD_{x_k} +\frac{1}{2}\sum_{j=1}^n \left(c_j(x)
h\partial_{x_j}+ h\partial_{x_j}\circ c_j(x)\right) + p_0(x) \\ \nonumber
& =:& P_2+iP_1+P_0,\quad D_{x_j} = \frac{1}{i}\frac{\partial}{\partial_{x_j}}.
\end{eqnarray}
Here the coefficients $b_{j,k}$, $c_j$, $p_0$ all belong to $C^{\infty}(\real^n;\real)$, with $b_{j,k}=b_{k,j}$.
Associated to $P$ in (\ref{eq1.1}) is the symbol in the semiclassical sense,
\begeq
\label{eq1.2}
p(x,\xi)=p_2(x,\xi)+ip_1(x,\xi)+p_0(x),
\endeq
\begeq
\label{eq1.3}
p_2(x,\xi)=\sum_{j,k=1}^n b_{j,k}(x) \xi_j \xi_k,\quad p_1(x,\xi)=\sum_{j=1}^n c_j(x)\xi_j,
\endeq
so that $p_j(x,\xi)$ is a real-valued polynomial in $\xi$, positively homogeneous of degree $j$, $0\leq j\leq 2$.
We may notice that $p(x,\xi)$ coincides with the Weyl symbol of $P$ modulo ${\cal O}(h^2)$, locally uniformly.

As in~\cite{HeHiSj}, we shall assume that
\begeq
\label{eq1.4}
p_2(x,\xi)\geq 0,\quad p_0(x)\geq 0.
\endeq
Furthermore, we shall impose the following growth conditions,
\begeq
\label{eq1.5}
\partial_x^{\alpha} b_{j,k}(x)={\cal O}(1),\quad \abs{\alpha}\geq 0,
\endeq
\begeq
\label{eq1.6}
\partial_x^{\alpha} c_j(x)={\cal O}(1),\quad \abs{\alpha}\geq 1,
\endeq
\begeq
\label{eq1.7}
\partial_x^{\alpha} p_0(x)={\cal O}(1),\quad \abs{\alpha}\geq 2.
\endeq
From section 3 in~\cite{HeHiSj}, we recall that under these assumptions, the graph closure of $P: {\cal S}(\real^n)\rightarrow
{\cal S}(\real^n)$, still denoted by $P$ and such that $\Re P\geq 0$, coincides with the maximal closed realization
of $P$, with the domain given by ${\cal D}(P)=\{u\in L^2;\, Pu\in L^2\}$. In particular, this shows that the operator $P$ is m-accretive, and hence,
the contraction semigroup
\begeq
\label{eq1.8}
e^{-tP/h}: L^2\rightarrow L^2,\quad t\geq 0
\endeq
is well-defined.

\medskip
We proceed next to recall the additional assumptions of a dynamical nature, introduced in section 4 of~\cite{HeHiSj}. Let
\begeq
\label{eq1.9}
\nu(x,\partial_x)=\sum_{j=1}^n c_j(x)\partial_{x_j},
\endeq
and recall the Hypothesis 4.1 of~\cite{HeHiSj}:
\begeq
\label{eq1.91}
\wrtext{The set}\,\{x\in \real^n;\, p_0(x)=0,\, \nu(x,\partial_x)=0\}\,\, \wrtext{is finite}\,\,
=\{x_1, \ldots\, x_N\}.
\endeq

\noindent
With $\rho_j=(x_j,0)$, $1\leq j\leq N$, we define the critical set
\begeq
\label{eq1.10}
{\cal C}=\{\rho_1,\ldots \, \rho_N\}\subset \real^{2n}.
\endeq
The coefficients $p_0$, $p_1$, $p_2$ in (\ref{eq1.2}) all vanish to the second order at each $\rho_j\in {\cal C}$.
As in~\cite{HeHiSj}, we define
\begeq
\label{eq1.11}
\widetilde{p}(x,\xi)=p_0(x)+\langle{\xi}\rangle^{-2} p_2(x,\xi),
\endeq
and consider the time average
\begeq
\label{eq1.12}
\langle{\widetilde{p}}\rangle_{T_0}=\frac{1}{T_0}\int_{-T_0/2}^{T_0/2} \widetilde{p}\circ \exp(tH_{p_1})\,dt,\quad T_0>0.
\endeq
We shall assume that the Hypothesis 4.3 of~\cite{HeHiSj} holds true:
\begeq
\label{eq1.12.5}
\wrtext{For}\,\,T_0>0\, \wrtext{fixed, we have near each}\, \rho_j,\,\,
\langle{\widetilde{p}\rangle}_{T_0}(\rho) \sim \abs{\rho-\rho_j}^2,
\endeq
and in any set of the form $\abs{x}\leq C$, ${\rm dist}(\rho,{\cal C})\geq 1/C$, we have
\begeq
\label{eq1.12.54}
\langle{\widetilde{p}\rangle}_{T_0}(\rho)\geq \frac{1}{\widetilde{C}(C)}
,\quad \widetilde{C}(C)>0.
\endeq
We also need an additional dynamical hypothesis near $\infty$ in
$\real^n$,
\begin{equation} \label{eq1.12.546}
\begin{split}
&\forall\hbox{ \neigh{} }U\hbox{ of }\pi _x{\cal C}, \hbox{ and }\forall \
x\in \real^n \sm U,\ \exists C>0, \\
&
{\rm meas\,}\left(\{ t\in [-{T_0\over 2},{T_0\over 2}];\, p_0(\exp t\nu
(x))\ge {1\over C}\}\right)\ge {1\over C}.
\end{split}
\end{equation}

\medskip
\noindent
Under the assumptions above, the paper~\cite{HeHiSj} defines an auxiliary real valued weight function
$\psi_{\eps}(x,\xi)={\cal O}(\eps)$ on $T^*\real^n$, $\eps>0$, such that
$$
\partial^{\alpha}_x\partial^{\beta}_{\xi} \psi_{\eps}(x,\xi)=
{\cal O}\left(\eps^{1-\abs{\alpha+\beta}/2}\langle{\xi}\rangle^{-\abs{\beta}}\right),
$$
together with an associated canonical transformation
\begeq
\label{eq1.12.55}
\kappa(\delta): \real^{2n}\rightarrow
\Lambda_{\delta}:=\{(x,\xi)+i\delta H_{\psi_{\eps}}(x,\xi);\, (x,\xi)\in \real^{2n}\},\quad 0<\delta \ll 1,
\endeq
for which
\begeq
\label{eq1.12.6}
\kappa(\delta)(x,\xi)=(x,\xi)+i\delta H_{\psi_{\eps}}(x,\xi)+{\cal O}(\eps^{1/2}\delta^2).
\endeq
Here we let $H_f$ denote the Hamilton field $f'_{\xi}(x,\xi)\cdot \partial_x - f'_{x}(x,\xi)\cdot \partial_{\xi}$ of
a $C^1$--function $f(x,\xi)$. We refer to section 4 of~\cite{HeHiSj} for the details of the construction of $\psi_{\eps}$ and $\kappa_{\delta}$.
Here we shall merely recall that $\delta>0$ fixed in (\ref{eq1.12.6}) should be small enough, and $\eps=Ah$, with $A$ arbitrarily large but fixed.

Associated to $\kappa(\delta)$ in (\ref{eq1.12.55}) there is an elliptic Fourier integral operator with a complex
phase $A_{\delta,\eps}$, constructed in Section 5 of~\cite{HeHiSj}, such that
\begeq
\label{eq1.12.7}
A_{\delta,\eps}: {\cal S}\rightarrow {\cal S}
\endeq
continuously, and
\begeq
\label{eq1.12.8}
A_{\delta,\eps}={\cal O}_A(1): L^2\rightarrow L^2.
\endeq
Moreover, it is proved in~\cite{HeHiSj} that $A_{\delta,\eps}$ is invertible, when  $\eps/h\gg 1$,
with the inverse $A_{\delta,\eps}^{-1}$ also enjoying the mapping properties (\ref{eq1.12.7}), (\ref{eq1.12.8}).

When $B\geq A$ fixed is to be chosen, and $\widetilde{A}\gg B$ is large enough, as in sections 6,7
of~\cite{HeHiSj}, we shall consider the conjugated operator
\begeq
\label{eq1.13}
P^{\delta,\widetilde{\eps}}=A_{\delta,\widetilde{\eps}}^{-1} P A_{\delta,\widetilde{\eps}},\quad
\widetilde{\eps}=\widetilde{A}h,
\endeq
acting on $L^2(\real^n)$. It was then proved in~\cite{HeHiSj} that the real part of the symbol of
$P^{\delta,\widetilde{\eps}}$ is $\geq \frac{\delta\widetilde{\eps}}{C}-C h$ outside of
${\cal C}+B(0,\sqrt{\widetilde{\eps}})$, $C>0$. In the set
${\cal C}+B(0,\sqrt{\widetilde{\eps}})$, the symbol of $P^{\delta,\widetilde{\eps}}$ is independent of $\widetilde{\eps}$
modulo ${\cal O}(\widetilde{\eps}\left(\frac{h}{\widetilde{\eps}}\right)^{\infty})$ and is of the form
\begeq
\label{eq1.13.1}
\widetilde{P}^{\delta} \sim p^{\delta}+h^2r_2+\ldots\,
\endeq
where
\begeq
\label{eq1.14}
\abs{p^{\delta}(\rho)}\sim {\rm dist}(\rho,{\cal C})^2, \quad \Re p^{\delta}(\rho) \sim {\rm dist}(\rho, {\cal C})^2.
\endeq
Here we write $B(0,r)=\{\rho\in \real^{2n}; \abs{\rho}<r\}$, $r>0$.

\medskip
As has also been recalled in section 8 in~\cite{HeHiSj}, in the set ${\cal C}+B(0,\sqrt{\widetilde{\eps}})$, we have
$$
P^{\delta,\widetilde{\eps}}-\widetilde{P}^{\delta}=\widetilde{{\cal O}}
\left(\widetilde{\eps}\left(\frac{h}{\widetilde{\eps}}\right)^{\infty}\right),
$$
while when away from ${\cal C}+B(0,\sqrt{\widetilde{\eps}})$, we shall use that
\begeq
\label{eq1.14.1}
P^{\delta,\widetilde{\eps}}-\widetilde{P}^{\delta}=\widetilde{{\cal O}}\left(\widetilde{\eps}+\rho^2\right).
\endeq
Here, following~\cite{HeHiSj}, we use the notation $f_{\eps}=\widetilde{\cal O}(\eps)$ to express that
$$
\partial^{\alpha}_x\partial_{\xi}^{\beta}f_{\eps}(x,\xi)=
{\cal O}\left(\eps^{1-\abs{\alpha+\beta}/2}\langle{\xi}\rangle^{-\abs{\beta}}\right),
$$
for arbitrary multi-indices $\alpha$, $\beta\in \nat^n$.

\bigskip
We shall study estimates for the real part of the quadratic form associated to
the operator $P^{\delta,\widetilde{\eps}}$. The starting point here is Proposition 7.1 of~\cite{HeHiSj}: let
$0\leq {k}_{\eps}=\widetilde{{\cal O}}(\eps)$ be equal to $\eps$ in ${\cal C}+B(0,\sqrt{\eps})$ and
have its support in ${\cal C}+B(0,\sqrt{2\eps})$. Let
$$
K_{\eps}={\rm Op}_h(k_{\eps})
$$
stand for the Weyl quantization of $k_{\eps}(x,h\xi)$. It is then established in Proposition 7.1 of~\cite{HeHiSj} that
\begeq
\label{eq2.1}
\Re \left((P^{\delta,\eps}+K_{\eps})u|u\right)\geq \left(\frac{\delta \eps}{C}-Ch\right)\norm{u}^2,\quad u\in {\cal S},
\endeq
when $\eps=Ah$, $C>0$ is independent of $\delta$, $A$, and $h$ is small enough depending on these 2 parameters.

Rather than working with the estimate (\ref{eq2.1}), we shall use that
\begeq
\label{eq2.1.0}
\Re \left((P^{\delta,\widetilde{\eps}}+K_{\eps})u|u\right)\geq \left(\frac{\delta \eps}{C}-Ch\right)\norm{u}^2,
\quad u\in {\cal S},
\endeq
which is proved in exactly the same way as in section 7 of~\cite{HeHiSj}. Here we recall that
$\eps=Ah$, $\widetilde{\eps}=\widetilde{A}h$, $\widetilde{A}\gg A$. In what follows we shall
use that the estimate (\ref{eq2.1.0}) holds also for $u\in {\cal D}(P^{\delta,\widetilde{\eps}})=
\{u\in L^2; P^{\delta,\widetilde{\eps}}u\in L^2\}=A_{\delta,\widetilde{\eps}}^{-1} {\cal D}(P)$.

\medskip
Let
$$
\Pi_B=\Pi_B^{\delta,\widetilde{\eps}}
$$
be the spectral projection of $P^{\delta,\widetilde{\eps}}$ associated with the spectrum of $P^{\delta,\widetilde{\eps}}$
in the open disc $D(0,Bh)$. From Theorem 8.3 in~\cite{HeHiSj} we recall that the spectrum of
$P^{\delta,\widetilde{\eps}}$ in $D(0,Bh)$ is discrete, and the eigenvalues are of the form
$$
\lambda_{j,k}(h)\sim h\left(\mu_{j,k}+h^{1/N_{j,k}}\mu_{j,k,1}+h^{2/N_{j,k}}\mu_{j,k,2}+\ldots\,\right),
$$
where the $\mu_{j,k}$ are all numbers in $D(0,B)$ of the form
\begeq
\label{eq2.1.1}
\mu_{j,k}=\frac{1}{i}\sum_{\ell=1}^n \left(\nu_{j,k,\ell}+\frac{1}{2}\right)\lambda_{j,\ell},\quad \nu_{j,k,\ell}\in \nat,
\endeq
for some $j\in \{1,\ldots N\}$. Here $\lambda_{j,\ell}$, $1\leq \ell \leq n$, are the eigenvalues of the Hamilton map of
the quadratic part of $p$ at $\rho_j\in {\cal C}$, for which
$\textrm{Im\,} \lambda_{j,\ell}>0$. Here we also assume that $B$ is chosen
such that $\abs{\mu_{j,k}}\neq B$, for all $j$, $k$.

Assume that $u\in L^2$ is such that
\begeq
\label{eq2.2}
u\in {\rm Ran}(1-\Pi_B).
\endeq
We are interested in lower bounds for
\begeq
\label{eq2.3}
\Re(P^{\delta,{\widetilde{\eps}}}u|u),
\endeq
which, in view of (\ref{eq2.1.0}), amounts to estimating $K_{\eps} u$. In doing so, we shall assume,
for notational simplicity only, that the critical set ${\cal C}$ defined in (\ref{eq1.10}) consists of a single point,
$\rho_1={(0,0)}$. From (\ref{eq1.14}), we know that the leading symbol
of $\widetilde{P}^{\delta}$, $p^{\delta}$, is such that
$$
\abs{p^{\delta}(\rho)}\sim \abs{\rho}^2,\quad \rho\in
B(0,\sqrt{\widetilde{\eps}}).
$$
Let
$$
p_0(x,\xi)=p_0^{\delta}(x,\xi)=\sum_{\abs{\alpha+\beta}=2}
\frac{\partial_x^{\alpha}\partial_{\xi}^{\beta}p^{\delta}(0,0)}{\alpha!\beta!} x^{\alpha} \xi^{\beta}
$$
be the quadratic approximation of $p^{\delta}$, so that
\begeq
\label{eq2.3.1}
p^{\delta}-p_0={\cal O}((x,\xi)^3)={\cal
O}\left((h+(x,\xi)^2)^{3/2}\right),\quad (x,\xi)\rightarrow (0,0).
\endeq
Then $p_0$ is an elliptic quadratic form on $\real^{2n}$, with a positive definite real part. The quadratic differential
operator
$$
P_0={\rm Op}_h(p_0): L^2\rightarrow L^2,
$$
has discrete spectrum, and from~\cite{Sj74} we know that the eigenvalues of $P_0$ are of the form $h\mu_{1,k}$,
with $\mu_{1,k}$ defined as in (\ref{eq2.1.1}).

\medskip
\noindent
When estimating $K_{\eps} u$ for $u\in {\rm Ran}(1-\Pi_B)$, we also introduce the spectral projection $\Pi_{0,B}$
associated to $P_0$ and the spectrum of $P_0$ in $D(0,Bh)$. Then, since $\Pi_Bu=0$,
\begeq
\label{eq2.31}
K_{\eps}u = K_{\eps}(\Pi_{0,B}-\Pi_B)u+K_{\eps}(1-\Pi_{0,B})u.
\endeq

We shall estimate the first term in the right hand side of (\ref{eq2.31}), using the following result.

\begin{lemma}
We have
\begeq
\label{eq2.32}
\Pi_B-\Pi_{0,B}={\cal O}_{B}(\widetilde{A}^{3/2}h^{1/2}+\widetilde{A}^{-1}): L^2\rightarrow L^2.
\endeq
\end{lemma}

\begin{proof}
Let $\chi \in C^{\infty}_0(B(0,2))$, $0 \leq \chi \leq 1$, be such
that $\chi(x,\xi)=1$ for $\abs{(x,\xi)}\leq 1$. Set
$\chi_{\sqrt{\widetilde{\eps}}}(x,\xi)=\chi(\widetilde{\eps}^{-1/2}(x,\xi))$. We shall first show that
\begeq
\label{eq2.35}
\Pi_{0,B}\left(1-\chi_{\sqrt{\widetilde{\eps}}}\right)={\cal
O}_B\left(\left(\frac{h}{\widetilde{\eps}}\right)^{\infty}\right): L^2\rightarrow
L^2,
\endeq
and similarly, that
\begeq
\label{eq2.4}
\Pi_B\left(1-\chi_{\sqrt{\widetilde{\eps}}}\right)={\cal
O}_B\left(\left(\frac{h}{\widetilde{\eps}}\right)^{\infty}\right):
L^2\rightarrow L^2.
\endeq
When proving (\ref{eq2.35}), we shall use the well-posed Grushin problem
for the quadratic operator $P_0$, described in~\cite{HeSjSt},~\cite{HeHiSj}. Let
\begeq
\label{eq2.41}
\Lambda=\biggl \langle{\frac{(x,hD)}{h^{1/2}}\biggr \rangle}=\left(1+\frac{x^2+(hD_x)^2}{h}\right)^{1/2},
\endeq
so that the quadratic elliptic operator $P_0$ is equipped with the natural
domain
$$
{\cal D}(P_0)=\{u\in L^2; \Lambda^2 u\in L^2\}.
$$
In section 11 in~\cite{HeSjSt}, using the generalized eigenfunctions of $P_0$ and of the adjoint $P_0^*$,
the authors have constructed the operators
\begeq
\label{eq2.42}
R_-:\comp^{N_0}\rightarrow L^2,\quad R_+: L^2\rightarrow \comp^{N_0},\quad N_0\in \nat,
\endeq
such that when $z\in D(0,Bh)$, the problem
\begeq
\label{eq2.42.1}
\left(P_0-z\right)u+R_-u_-=v,\quad R_+u=v_+,
\endeq
for $v\in L^2$, $v_+\in \comp^{N_0}$ has a unique solution $u\in {\cal D}(P_0)$, $u_-\in \comp^{N_0}$. Moreover, we have
the a priori estimate
\begeq
\label{eq2.42.2}
h\norm{\Lambda^2 u}+\abs{u_-}\leq {\cal O}(1)\left(\norm{v}+h\abs{v_+}\right).
\endeq

\medskip
Associated to (\ref{eq2.42.1}) is the Grushin operator
\begeq
\label{eq2.43}
{\cal P}_0(z)=
\left( \begin{array}{ccc}
P_0-z & R_- \\\
R_+ & 0
\end{array} \right): {\cal D}(P) \times \comp^{N_0} \rightarrow
L^2 \times \comp^{N_0},\quad z\in D(0,Bh),
\endeq
with an inverse
\begeq
\label{eq2.44}
{\cal E}_0(z)=\left( \begin{array}{ccc}
E(z) & E_+(z) \\\
E_-(z) & E_{-+}(z)
\end{array} \right): L^2 \times \comp^{N_0} \rightarrow {\cal D}(P) \times \comp^{N_0},
\endeq
depending holomorphically on $z$. From section 11 of~\cite{HeSjSt}, we recall that ${\cal E}_0(z)$
enjoys the following localization properties, when $k\in \real$,
\begeq
\label{eq2.44.1}
\Lambda^{2-k} E(z) \Lambda^k = {\cal O}\left(\frac{1}{h}\right): L^2\rightarrow L^2,
\endeq
and
\begeq
\label{eq2.45}
\Lambda^kE_+(z)={\cal O}_k(1): \comp^{N_0}\rightarrow L^2,\quad
E_-(z) \Lambda^{k}={\cal O}_k(1): L^2\rightarrow \comp^{N_0},
\endeq

Let $\gamma\subset D(0,B)$ be a simple positively oriented closed $h$-independent contour, such that all eigenvalues of $P_0$ and
$P^{\delta,\widetilde{\eps}}$ in $D(0,Bh)$ are contained in the interior of $h\gamma$, so that we have
$$
{\rm dist}(z,{\rm Spec}(P_0)\cup {\rm Spec}(P^{\delta,\widetilde{\eps}}))\geq
h/{\cal O}(1),\quad z\in h\gamma.
$$
Here we continue to assume that $B>0$ is chosen so that there are no numbers of the form $\mu_{j,k}$ in (\ref{eq2.1.1})
on the boundary of $D(0,B)$. Writing
\begeq
\label{eq2.46}
\Pi_{0,B}=\frac{1}{2\pi i} \int_{h \gamma} (z-P_0)^{-1}\,dz
\endeq
and using the well-known formula
\begeq
\label{eq2.5}
\left(z-P_0\right)^{-1}=-E(z)+E_+(z)E_{-+}(z)^{-1} E_-(z),
\endeq
we obtain that
\begeq
\label{eq2.6}
\Pi_{0,B}=\frac{1}{2\pi
i}\int_{h \gamma} E_+(z) E_{-+}(z)^{-1} E_-(z)\,dz.
\endeq
Now (\ref{eq2.45}) gives that for each $k \in \nat$,
\begeq
\label{eq2.8}
E_-(z)(1-\chi_{\sqrt{\widetilde{\eps}}})=E_-\Lambda^k \Lambda^{-k}
(1-\chi_{\sqrt{\widetilde{\eps}}})={\cal
O}\left(\left(\frac{h}{\widetilde{\eps}}\right)^{k/2}\right): L^2\rightarrow
\comp^{N_0}.
\endeq
Using also that along $h\gamma$, we have
$$
E_+(z)={\cal O}(1): \comp^{N_0}\rightarrow L^2,
$$
and
$$
E_{-+}(z)^{-1}={\cal O}(h^{-1}): \comp^{N_0}\rightarrow \comp^{N_0},
$$
as well as the fact that length of $h \gamma$ is ${\cal O}_B(h)$, we obtain (\ref{eq2.35}).

\medskip
\noindent
The proof of (\ref{eq2.4}) proceeds along the similar lines, relying upon the well-posed Grushin problem for
$P^{\delta,\widetilde{\eps}}-z$, constructed from the Grushin problem for $P_0$ and described in detail in section 11
of~\cite{HeSjSt} and section 8 of~\cite{HeHiSj}. In particular, the analogue of the localization property (\ref{eq2.45})
holds true for the inverse of the Grushin operator for $P^{\delta,\widetilde{\eps}}$, and arguing as above, we get
(\ref{eq2.4}).

\medskip
\noindent
We now come to consider estimates for the difference $\left(\Pi_B-\Pi_{0,B}\right)\chi_{\sqrt{\widetilde{\eps}}}$, where
we claim that
\begeq
\label{eq2.85}
\left(\Pi_B-\Pi_{0,B}\right)\chi_{\sqrt{\widetilde{\eps}}}={\cal
O}_B\left(\frac{\widetilde{\eps}^{3/2}}{h}+\frac{1}{\widetilde{A}}\right): L^2\rightarrow L^2,\quad \widetilde{\eps}=\widetilde{A}h.
\endeq
In view of (\ref{eq2.35}) and (\ref{eq2.4}), this will complete the proof of Lemma 2.1.

\medskip
\noindent
When proving (\ref{eq2.85}), we shall first estimate the difference
$\left(\widetilde{\Pi}_B-\Pi_{0,B}\right)\chi_{\sqrt{\widetilde{\eps}}}$,
where $\widetilde{\Pi}_B$ is the spectral projection
of the operator $\widetilde{P}^{\delta}$ associated with the spectrum of $\widetilde{P}^{\delta}$ in $D(0,Bh)$.
Using \ref{eq2.46}) together with the similar formula for $\widetilde{\Pi}_B$, we get, by an application of
the resolvent identity,
\begeq
\label{eq2.9}
\left(\widetilde{\Pi}_B-\Pi_{0,B}\right)\chi_{\sqrt{\widetilde{\eps}}}=\frac{1}{2\pi i}\int_{h \gamma}
\left(z-\widetilde{P}^{\delta}\right)^{-1}
\left(\widetilde{P}^{\delta}-P_0\right)\left(z-P_0\right)^{-1}\chi_{\sqrt{\widetilde{\eps}}}\,dz.
\endeq
Now (\ref{eq2.44.1}), (\ref{eq2.45}), (\ref{eq2.5}), together with the fact that along $h\gamma$ we have
$\norm{E_{-+}^{-1}(z)}={\cal O}(h^{-1})$, imply that for $k\in \real$,
\begeq
\label{eq2.10}
\Lambda^k \left(z-P_0\right)^{-1} \Lambda^{-k} = {\cal
O}\left(\frac{1}{h}\right) :L^2\rightarrow L^2,
\endeq
and even that
\begeq
\label{eq2.10.1}
\Lambda^{2+k} \left(z-P_0\right)^{-1} \Lambda^{-k} = {\cal O}\left(\frac{1}{h}\right): L^2 \rightarrow L^2.
\endeq
Here we shall take $k=3$ in (\ref{eq2.10}). Writing the integrand in (\ref{eq2.9}) as
\begin{eqnarray}
\label{eq2.11}
& & \left(z-\widetilde{P}^{\delta}\right)^{-1}\left(\widetilde{P}^{\delta}-P_0\right) \left(z-P_0\right)^{-1}
\chi_{\sqrt{\widetilde{\eps}}} \\ \nonumber
& = & \left(z-\widetilde{P}^{\delta}\right)^{-1}\left(\widetilde{P}^{\delta}-P_0\right)\Lambda^{-3} \Lambda^3
(z-P_0)^{-1} \Lambda^{-3} \Lambda^3 \chi_{\sqrt{\widetilde{\eps}}}.
\end{eqnarray}
we see that we have to estimate the operator norm of $(\widetilde{P}^{\delta}-P_0)\Lambda^{-3}$. Now it follows
from (\ref{eq1.13.1}), (\ref{eq2.3.1}) that
\begeq
\label{eq2.12}
(\widetilde{P}^{\delta}-P_0)\Lambda^{-3}={\cal O}(h^{3/2}): L^2\rightarrow L^2.
\endeq
Also,
\begeq
\label{eq2.13}
\Lambda^3 \chi_{\sqrt{\widetilde{\eps}}}={\cal
O}\left(\frac{\widetilde{\eps}^{3/2}}{h^{3/2}}\right): L^2\rightarrow L^2.
\endeq
Combining (\ref{eq2.9}), (\ref{eq2.10}), (\ref{eq2.12}), (\ref{eq2.13})
together with the fact that
$$
(z-\widetilde{P}^{\delta})^{-1}={\cal O}(h^{-1}): L^2\rightarrow L^2,\quad z\in
h\gamma,
$$
which follows from Theorem 8.4 in~\cite{HeHiSj},
and that the length of $h\gamma$ is ${\cal O}_B(h)$, we obtain that
\begeq
\label{eq2.14}
\left(\widetilde{\Pi}_B-\Pi_{0,B}\right)\chi_{\sqrt{\widetilde{\eps}}}=
{\cal O}_B\left(\frac{\widetilde{\eps}^{3/2}}{h}\right): L^2\rightarrow L^2.
\endeq

\medskip
It only remains now to estimate the operator norm of $\left(\Pi_B-\widetilde{\Pi}_B\right)
\chi_{\sqrt{\widetilde{\eps}}}$. To that end, we write, as in
(\ref{eq2.9}),
\begin{eqnarray}
\label{eq2.15}
\left(\Pi_B-\widetilde{\Pi}_B\right)\chi_{\sqrt{\widetilde{\eps}}} & =
& \frac{1}{2\pi
i} \int_{h\gamma} \left(z-P^{\delta,\widetilde{\eps}}\right)^{-1}
\left(P^{\delta,\widetilde{\eps}}-\widetilde{P}^{\delta}\right)\left(z-\widetilde{P}^{\delta}\right)^{-1}
\chi_{\sqrt{\widetilde{\eps}}}\,dz \\ \nonumber
& = & \frac{1}{2\pi i} \int_{h\gamma} \left(z-P^{\delta,\widetilde{\eps}}\right)^{-1}
\left(P^{\delta,\widetilde{\eps}}-\widetilde{P}^{\delta}\right)\chi_{\sqrt{\widetilde{\eps}}}
\left(z-\widetilde{P}^{\delta}\right)^{-1}
\chi_{\sqrt{\widetilde{\eps}}}\,dz \\ \nonumber
& & + \frac{1}{2\pi i} \int_{h\gamma} \left(z-P^{\delta,\widetilde{\eps}}\right)^{-1}
\left(P^{\delta,\widetilde{\eps}}-\widetilde{P}^{\delta}\right)\left(1-\chi_{\sqrt{\widetilde{\eps}}}\right)
\left(z-\widetilde{P}^{\delta}\right)^{-1}
\chi_{\sqrt{\widetilde{\eps}}}\,dz\\ \nonumber
& = & {\rm I}+{\rm II},
\end{eqnarray}
with the natural definitions of ${\rm I}$ and ${\rm II}$. Using, as in section 8 of~\cite{HeHiSj}, that
$$
\left(P^{\delta,\widetilde{\eps}}-\widetilde{P}^{\delta}\right)\chi_{\sqrt{\widetilde{\eps}}} =
{\cal O}\left(h\frac{h}{\widetilde{\eps}}\right): L^2\rightarrow L^2,
$$
together with the ${\cal O}(h^{-1})$--estimates for the resolvents of
$P^{\delta,\widetilde{\eps}}$ and $\widetilde{P}^{\delta}$ along the
contour $h\gamma$, we get
\begeq
\label{eq2.16}
{\rm I} = {\cal O}_B\left(\frac{1}{h}{\widetilde{\eps}}\right)={\cal O}_B\left(\frac{1}{\widetilde{A}}\right):
L^2\rightarrow L^2.
\endeq

\medskip
\noindent
We now come to estimate the term ${\rm II}$ in (\ref{eq2.15}). We have
\begeq
\label{eq2.16.01}
\chi_{\sqrt{\widetilde{\eps}}}\left(P^{\delta,\widetilde{\eps}}-\widetilde{P}^{\delta}\right)={\cal
O}\left(h\frac{h}{\widetilde{\eps}}\right): L^2\rightarrow L^2,
\endeq
and using this estimate, as well as the ${\cal O}(h^{-1})$--resolvent bounds for $P^{\delta,\widetilde{\eps}}$ and
$\widetilde{P}^{\delta}$, we see that modulo a term whose operator norm on $L^2$ is
$$
{\cal O}_B\left(\frac{h}{\widetilde{\eps}}\right)={\cal O}_B\left(\widetilde{A}^{-1}\right),
$$
we may replace the integrand in ${\rm II}$ by the following expression,
\begeq
\label{eq2.16.02}
\left(z-P^{\delta,\widetilde{\eps}}\right)^{-1}
\left(1-\chi_{\sqrt{\widetilde{\eps}}}\right)\left(P^{\delta,\widetilde{\eps}}-\widetilde{P}^{\delta}\right)
\left(1-\chi_{\sqrt{\widetilde{\eps}}}\right)\left(z-\widetilde{P}^{\delta}\right)^{-1}
\chi_{\sqrt{\widetilde{\eps}}}.
\endeq
Here following (\ref{eq2.11}), we write
\begeq
\label{eq2.16.03}
\left(P^{\delta,\widetilde{\eps}}-\widetilde{P}^{\delta}\right)
\left(1-\chi_{\sqrt{\widetilde{\eps}}}\right)\left(z-\widetilde{P}^{\delta}\right)^{-1}
\chi_{\sqrt{\widetilde{\eps}}}
\endeq
as
\begeq
\label{eq2.16.04}
\left(P^{\delta,\widetilde{\eps}}-\widetilde{P}^{\delta}\right)
\left(1-\chi_{\sqrt{\widetilde{\eps}}}\right)\Lambda^{-k-2}
\left(\Lambda^{k+2}
\left(z-\widetilde{P}^{\delta}\right)^{-1}\Lambda^{-k}\right)\Lambda^k
\chi_{\sqrt{\widetilde{\eps}}}.
\endeq
Here using (\ref{eq1.14.1}) we see that
\begeq
\label{eq2.16.06}
\left(P^{\delta,\widetilde{\eps}}-\widetilde{P}^{\delta}\right)
\left(1-\chi_{\sqrt{\widetilde{\eps}}}\right)\Lambda^{-k-2} = {\cal
O}\left(\frac{h^{k/2+1}}{\widetilde{\eps}^{k/2}}\right).
\endeq
On the other hand, as in (\ref{eq2.10.1}),
\begeq
\label{eq2.16.07}
\Lambda^{k+2}\left(z-\widetilde{P}^{\delta}\right)^{-1}\Lambda^{-k} =
{\cal O}\left(\frac{1}{h}\right): L^2 \rightarrow L^2,
\endeq
and also,
\begeq
\label{eq2.16.08}
\Lambda^k \chi_{\sqrt{\widetilde{\eps}}} = {\cal
O}\left(\frac{\widetilde{\eps}^{k/2}}{h^{k/2}}\right): L^2\rightarrow L^2.
\endeq
Combining (\ref{eq2.16.04}), (\ref{eq2.16.06}), (\ref{eq2.16.07}), and (\ref{eq2.16.08}), we see that the
expression (\ref{eq2.16.03}) is ${\cal O}(1)$.

\bigskip
\noindent
We now come to estimate the remaining factor in (\ref{eq2.16.02}). To that end, we let
$$
L_{\widetilde{\eps}} = \widetilde{{\cal O}}\left(1+\frac{{\rm min}((x,\xi)^2,\widetilde{\eps})}{h}\right)
$$
be an elliptic symbol in the class defined by the right hand side, and write
\begeq
\label{eq2.16.09}
\left(z-P^{\delta,\widetilde{\eps}}\right)^{-1}
\left(1-\chi_{\sqrt{\widetilde{\eps}}}\right)=\left(
\left(z-P^{\delta,\widetilde{\eps}}\right)^{-1} L_{\widetilde{\eps}} \right)
L_{\widetilde{\eps}}^{-1}\left(1-\chi_{\sqrt{\widetilde{\eps}}}\right).
\endeq
Here we know that
$$
\left(z-P^{\delta,\widetilde{\eps}}\right)^{-1} L_{\widetilde{\eps}} = {\cal
O}\left(\frac{1}{h}\right): L^2\rightarrow L^2,
$$
and since
$$
L_{\widetilde{\eps}}^{-1} \left(1-\chi_{\sqrt{\widetilde{\eps}}}\right) = {\cal
O}\left(\frac{h}{\widetilde{\eps}}\right): L^2 \rightarrow L^2,
$$
it follows that the expression (\ref{eq2.16.09}) is ${\cal
O}\left(\frac{1}{\widetilde{\eps}}\right)$. Using finally that the length
of the integration contour in (\ref{eq2.15}) is ${\cal O}_B(h)$, we get
\begeq
\label{eq2.17}
{\rm II} = {\cal O}_B\left(\frac{h}{\widetilde{\eps}}\right)={\cal
O}_B\left(\frac{1}{\widetilde{A}}\right): L^2\rightarrow L^2.
\endeq

\medskip
\noindent
Combining (\ref{eq2.15}), (\ref{eq2.16}), (\ref{eq2.17}), we conclude that
$$
\left(\Pi_B-\widetilde{\Pi}_B\right)\chi_{\sqrt{\widetilde{\eps}}}={\cal
O}_B\left(\frac{1}{\widetilde{A}}\right): L^2\rightarrow L^2.
$$
In view of (\ref{eq2.14}), the bound (\ref{eq2.85}) follows, and this completes the proof of Lemma 2.1.
\end{proof}

\bigskip
We now come to estimate the second term in the right hand side of
(\ref{eq2.31}), given by $K_{\eps}(1-\Pi_{0,B})u$. Here, the difficulty is that in general, due to a pseudospectral
phenomenon~\cite{DeSjZw}, the operator norm of $\Pi_{0,B}$ may exhibit some exponential growth, as
$h\rightarrow 0$. To circumvent this issue, our fist task will be to establish a more manageable characterization of
the vector $v=(1-\Pi_{0,B})u$. Specifically, we shall now discuss properties of the range of the projection
$1-\Pi_{0,B}$ on $L^2$.

\medskip
\noindent
In (\ref{eq2.1.1}), following~\cite{Sj74}, we have already recalled the form of the eigenvalues
of the elliptic quadratic operator $P_0$. From~\cite{Sj74}, we know furthermore
that the generalized eigenfunctions of $P_0$ are of the form
\begeq
\label{eq2.18}
h^{-n/4} p\left(\frac{x}{\sqrt{h}}\right)e^{i\Phi(x)/h},
\endeq
where $p(x)$ is a polynomial on $\real^n$ and $\Phi(x)$ is a quadratic
form with $\textrm{Im\,} \Phi>0$. The degree of the polynomial $p(x)$ in
(\ref{eq2.18}) tends to $\infty$ together with the real part of the
eigenvalue $h \mu_{j,k}$ in (\ref{eq2.1.1}). We may also recall from~\cite{Sj74}
that $\Phi$ in (\ref{eq2.18}) is such that the positive Lagrangian subspace $\Lambda_{\Phi}=\{(x,\Phi'(x)),
x\in \comp^n\}$ is the direct sum of the generalized eigenspaces of
the Hamilton map of $p_0$, corresponding to the eigenvalues with
a positive imaginary part. Correspondingly, the generalized
eigenfunctions of the formal $L^2$ adjoint $P_0^*$ are of the form
\begeq
\label{eq2.19}
h^{-n/4} q\left(\frac{x}{\sqrt{h}}\right)e^{i\Psi(x)/h},
\endeq
where $q$ is a polynomial and $\Psi$ is a quadratic form such that $\textrm{Im\,} \Psi$ is
positive definite.

Let $e_1,\ldots e_N$ be a basis for ${\rm Ran}\,(\Pi_{0,B})$ and let $e_1^*,\ldots e_N^*$ be the
corresponding dual basis for ${\rm Ran}\, \left((\Pi_{0,B})^*\right)$. If $v\in L^2$, we have
\begeq
\label{eq2.20}
\Pi_{0,B} v = \sum_{j=1}^N (v|e_j^*)e_j,
\endeq
and therefore, $v\in {\rm Ran}(1-\Pi_{0,B})$ precisely when
$v$ is orthogonal to ${\rm Ran} \left((\Pi_{0,B})^*\right)$.

\begin{prop}
There exists a selfadjoint $h$--differential operator $Q={\rm
Op}_h(q)$, where $q$ is a positive definite quadratic form on $T^*\real^n$, such that
\begeq
\label{eq2.21}
{\rm Ran} \left(E\left(Q,\frac{Bh}{C}\right)\right)\subset {\rm Ran}((\Pi_{0,B})^*)\subset {\rm
Ran}(E(Q,CBh)),
\endeq
for some $C>1$ which is independent of $Q$ and $B$. Here
$E(Q,\lambda)=1_{(-\infty,\lambda]}(Q)$ is the finite rank spectral projection
associated to $Q$ and the interval $(-\infty,\lambda]$.
\end{prop}

\begin{proof}
The operator $Q$ will be seen to be essentially the $h$--Weyl
quantization of the classical harmonic oscillator on $\real^n$.
When constructing $Q$, recall that the generalized eigenfunctions of
$P_0^*$ are of the form (\ref{eq2.19}). We shall write $\Psi(x)=(Bx,x)$, where $B$
is a symmetric matrix, $B=B_1+iB_2$, where $B_j$ are real, $j=1,2$,
and $B_2>0$. The real linear canonical transformation
\begeq
\label{eq2.22}
\kappa_1: (x,\xi)\mapsto (x,\xi-B_1 x)
\endeq
maps the positive Lagrangian subspace $\Lambda_{\Psi}=\{(x,Bx); x\in \comp^n\}$ to the
positive Lagrangian subspace $\{(x,iB_2 x); x\in \comp^n\}$. Now since $B_2>0$, there exists an invertible
real $n\times n$ matrix $C$ such that the real linear canonical transformation
\begeq
\label{eq2.23}
\kappa_2: (x,\xi)\mapsto (C^{-1}x, C^t \xi)
\endeq
maps $\{(x,iB_2x); x\in \comp^n\}=\kappa_1(\Lambda_{\Psi})$ to $\{(x,ix); x\in \comp^n\}$. We take the operator
\begeq
\label{eq2.24}
\widetilde{Q}=\frac{1}{2}\sum_{j=1}^n
\left(x_j^2+(hD_{x_j})^2\right)={\rm Op}_h(\widetilde{q}),\quad
\widetilde{q}(x,\xi)=\frac{1}{2}\sum_{j=1}^n \left(x_j^2+\xi_j^2\right),
\endeq
associated to $\Lambda_{\widetilde{\varphi}}=\{(x,ix); x\in \comp^n\}=\left(\kappa_2\circ \kappa_1\right)\left(\Lambda_{\Psi}\right)$,
$\widetilde{\varphi}(x)=ix^2/2$. To obtain the operator $Q$ it only
remains to notice that associated to $\kappa_1$ and $\kappa_2$ we have
the metaplectic operators
\begeq
\label{eq2.25}
U_1: L^2\rightarrow L^2,\quad U_1 f(x)=e^{-i(B_1x,x)/2h} f(x),
\endeq
and
\begeq
\label{eq2.26}
U_2: L^2\rightarrow L^2,\quad U_2 f(x)=f(Cx) \abs{{\rm det}\, C}^{1/2},
\endeq
both unitary on $L^2$, and hence with $U:=U_2\circ U_1$, we can take
$Q:=U^{-1} \widetilde{Q} U={\rm Op}_h\left(\left(\widetilde{q}\circ \left(\kappa_2\circ
\kappa_1\right)\right)\right)$. Notice that the eigenfunctions of $Q$
are of the form
\begeq
\label{eq2.26.2}
e_{\alpha,h}(x)=e_{\alpha,h=1}\left(\frac{x}{\sqrt{h}}\right),\quad
e_{\alpha,h=1}(x)=H_{\alpha}(C^{-1}x) e^{i\Psi(x)},\quad Q e_{\alpha,h}=h\left(\abs{\alpha}+\frac{n}{2}\right) e_{\alpha,h},
\endeq
where $H_{\alpha}(x)=\prod_{j=1}^n H_{\alpha_j}(x_j)$, $\alpha \in
\nat^n$, are the Hermite polynomials. The result follows.
\end{proof}

Having established a favorable comparison for the linear space
${\rm Ran} \left(\left(\Pi_{0,B}\right)^*\right)\subset L^2$, we return to the problem
of estimating $K_{\eps}v$, for $v=(1-\Pi_{0,B})u$. Let $\psi\in
C^{\infty}(\real; [0,1])$ with ${\rm supp}(\psi)\subset [0,1]$, and
set
\begeq
\label{eq2.26.5}
\psi_{\lambda}(t)=\psi\left(\frac{t}{\lambda}\right),\quad \lambda>0.
\endeq
It follows then from Proposition 2.2 that
\begeq
\label{eq2.27}
\psi_{\frac{Bh}{C}}(Q)v=0.
\endeq
To understand the operator occurring in (\ref{eq2.27}), it is convenient to perform a suitable dilation in
phase space. Assume therefore that $\lambda>0$ in (\ref{eq2.26.5}) is such that $h\ll \lambda \ll 1$.
Let us make the change of variables
$$
x=\lambda^{1/2}\widetilde{x},\quad D_{x}=\lambda^{-1/2}D_{\widetilde{x}}.
$$
Then, since the operator $Q$ is quadratic,
\begeq
\label{eq2.27.1}
\frac{1}{\lambda}Q=\frac{1}{\lambda} q^w(x,hD_x)=\frac{1}{\lambda}
q^w\left(\lambda^{1/2}\left(\widetilde{x},\frac{h}{\lambda}D_{\widetilde{x}}\right)\right)=
q^w\left(\widetilde{x},\frac{h}{\lambda}D_{\widetilde{x}}\right).
\endeq
It follows therefore from the functional calculus in the version of~\cite{DiSj} that
\begeq
\label{eq2.27.2}
\psi(\lambda^{-1}Q)={\rm Op}_{\frac{h}{\lambda},\widetilde{x}}
\left(r\left(\widetilde{x},\widetilde{\xi}; \frac{h}{\lambda}\right)\right)=
r^w\left(\lambda^{-1/2}\left(x,hD_x\right);\frac{h}{\lambda}\right),
\endeq
where $r\in S(\langle{\cdot \rangle}^{-N})$ for any $N\in \nat$, with
a complete asymptotic expansion in each of these symbol spaces,
and with the leading symbol $\psi(q(\widetilde{x},\widetilde{\xi}))$.

\medskip
\noindent
{\it Remark}. It is well known~\cite{Ho} that when $Q={\rm Op}_h^w(q)$
where $q$ is a positive definite quadratic form, then the
Weyl symbol of $f(Q)$, $f\in C^{\infty}_0(\real)$, is of the form $\widetilde{f}(q;h)=f(q)+{\cal O}(h)$. It follows
therefore that in (\ref{eq2.27.2}) we have
$$
r\left(\widetilde{x},\widetilde{\xi};\frac{h}{\lambda}\right) =
\widetilde{\psi}\left(q(\widetilde{x},\widetilde{\xi});\frac{h}{\lambda}\right),
$$
where the leading term of $\widetilde{\chi}$ is
$\psi(q(\widetilde{x},\widetilde{\xi}))$.

\bigskip
It is therefore clear that in (\ref{eq2.1}) we can take
\begeq
\label{eq2.28}
K_{\eps}=\eps \psi_{\frac{Bh}{C}}(Q),\quad \eps=Ah,
\endeq
for a suitable choice of $B\geq A$ fixed, where we can take $B$ as a
fixed multiple of $A$. With this choice, we get, using (\ref{eq2.27}),
\begeq
\label{eq2.29}
K_{\eps}(1-\Pi_{0,B})u=0.
\endeq

\bigskip
Combining (\ref{eq2.1.0}), (\ref{eq2.31}), Lemma 2.1, and (\ref{eq2.29}), we see that
for $u\in L^2$ with $u\in {\rm Ran}(1-\Pi_B)$, we have
\begeq
\label{eq24}
\left(\frac{\delta
\eps}{C}-Ch\right)\norm{u}^2 \leq \Re (P^{\delta,\widetilde{\eps}}u,u)+
\eps {\cal
O}_B\left(\widetilde{A}^{3/2}h^{1/2}+\widetilde{A}^{-1}\right)\norm{u}^2.
\endeq
Recall that here $B\geq A$, $B={\cal O}(A)$, is taken fixed, and $\delta>0$ is sufficiently small but fixed.
Choosing first $\widetilde{A}\gg B$ large enough and then taking $h$ sufficiently small depending on these parameters,
we absorb the second term in the right hand side of (\ref{eq24}) into the left hand side.

\begin{prop}
When $A\leq B\ll \widetilde{A}$, let $\Pi_B$ be the spectral projection of $P^{\delta,\widetilde{\eps}}$,
$\widetilde{\eps}=\widetilde{A}h$, associated with $D(0,Bh)$. Here $B$ is a fixed multiple of $A$.
Assume that $u\in {\cal D}(P^{\delta,\widetilde{\eps}})$ is such that
$u\in {\rm Ran}(1-\Pi_{B})$. Then for $h$ sufficiently small, we have
\begeq
\label{eq25}
\Re \left(P^{\delta,\widetilde{\eps}}u|u\right)\geq \frac{Bh}{{\cal
O}(1)}\norm{u}^2,\quad \widetilde{\eps}=\widetilde{A}h.
\endeq
\end{prop}

Now recall that
$$
P^{\delta,\widetilde{\eps}}=A_{\delta,\widetilde{\eps}}^{-1} P A_{\delta,\widetilde{\eps}},
$$
where $A_{\delta,\widetilde{\eps}}$,
$A_{\delta,\widetilde{\eps}}^{-1}: {\cal S}\rightarrow {\cal S}$, $L^2\rightarrow L^2$, have $L^2$ norm
${\cal O}_{\widetilde{A}}(1)$.
It is therefore clear from Proposition 2.3 that if $u$ is such that $u\in {\rm Ran}(1-\Pi)$, where
\begeq
\label{eq26}
\Pi=\frac{1}{2\pi i}\int_{h \gamma} \left(z-P\right)^{-1}\,dz
\endeq
is the spectral projection of $P$ associated with the spectrum of $P$ in $D(0,Bh)$, then
\begeq
\label{eq27}
\norm{e^{-tP/h}u}\leq {\cal O}(1)e^{-t/C}\norm{u},\quad C=C(B)>0.
\endeq
Therefore, it only remains to consider the restriction of the semigroup $e^{-tP/h}$ to the
finite-dimensional subspace ${\rm Ran}(\Pi)$, generated by the generalized eigenfunctions of $P$ corresponding to the
eigenvalues of $P$ of modulus $<Bh$. We shall now proceed to do so, in the framework of supersymmetric differential
operators.

\section{Supersymmetric operators and return to equilibrium in the double well case}
\setcounter{equation}{0}
The purpose of this section is to establish Theorem 1.1 in its general form, for a class of supersymmetric second order
differential operators, including (\ref{eq01}). Specifically, let
\begeq
\label{eq3.0}
A: \real^n\rightarrow \real^n
\endeq
be an invertible constant matrix. We decompose
\begeq
\label{eq3.001}
A=B+C,\quad ^tB=B,\quad ^tC=-C,
\endeq
and assume that
\begeq
\label{eq3.002}
B\geq 0.
\endeq
When $\phi\in C^{\infty}(\real^n;\real)$ is a Morse function such that
\begeq
\label{eq3.01}
\partial_x^{\alpha}\phi(x)={\cal O}(1),\quad
\partial_x^{\alpha}\left(\langle{B\partial_x\phi,\partial_x\phi\rangle}\right)=
{\cal O}(1),\quad \abs{\alpha}\geq 2,
\endeq
we consider the Witten-Hodge Laplacian associated to $A$ and $\phi$, acting on scalar functions,
defined as in section 10 of~\cite{HeHiSj},
\begin{eqnarray}
\label{eq3.02}
P = -\Delta_A^{(0)}=& & \sum_{j,k}hD_{x_j}B_{j,k} hD_{x_k}+
\sum_{j,k}\left(\partial_{x_j}\phi\right)B_{j,k}\left(\partial_{x_k}\phi\right)-h{\rm tr}(B\phi'') \\ \nonumber
& + & \sum_{j,k}\left(\left(\partial_{x_k}\phi\right)C_{j,k}h\partial_{x_j}+h\partial_{x_j}\circ
C_{j,k}\left(\partial_{x_k}\phi\right)\right).
\end{eqnarray}
The principal symbol of $P$ is of the form
\begeq
\label{eq3.1}
p(x,\xi)=\langle{B\xi,\xi}\rangle+2i\langle{C\phi'_x,\xi}\rangle+\langle{B\phi'_x,\phi'_x}\rangle,
\endeq
so that the assumptions (\ref{eq1.5}), (\ref{eq1.6}), (\ref{eq1.7}) are satisfied.

Assume that the Morse function $\phi$ has finitely many critical points $x_1,\ldots
x_N\in \real^n$ and that
\begeq \label{morse}
\abs{\phi'(x)}\geq \frac{1}{C},\quad \abs{x}\geq C.
\endeq
The assumption (\ref{eq1.91}) holds with
$$
{\cal C}=\{\rho_j; j=1,\ldots N\},\quad \rho_j=(x_j,0),
$$
and we shall also assume that the dynamical assumptions
(\ref{eq1.12.5}), (\ref{eq1.12.54}), and (\ref{eq1.12.546}) are valid. We then
know that the results of section 2 can be applied to $P$ in (\ref{eq3.02}).

\medskip
As in~\cite{HeHiSj}, we shall assume now that
\begin{eqnarray}
\label{eq3.1.01}
& & \phi\,\, \wrtext{has precisely three critical points, of which} \\ \nonumber
& & \wrtext{two are local minima}\,\,U_{\pm 1}, \wrtext{and the third one}\,\, U_0 \,\, \wrtext{is of index one.}
\end{eqnarray}
Then we know from~\cite{HeHiSj} that for $C>0$ large enough,
$P$ in (\ref{eq3.02}) has precisely 2 eigenvalues $\mu_0=0$
and $\mu_1$ in the disc $D(0,h/C)$ for $h$ small enough. Here $\mu_1$ is real and such that
\begeq
\label{eq3.1.1}
\mu_1=h\left(a_1(h) e^{-2S_{1}/h}+a_{-1}(h) e^{-2S_{-1}/h}\right), \quad S_j=\phi(U_0)-\phi(U_j)>0,\quad j=\pm 1,
\endeq
where $a_j(h)$ are real, $a_j(h)\sim a_{j,0}+h a_{j,1}+\ldots$, $a_{j,0}>0$.

The set $\phi^{-1}((-\infty,\phi(U_0)))$ has precisely two
connected components $D_j$, $j=\pm 1$, determined by the condition $U_j\in D_j$. Let $0\leq \chi_j\in C^{\infty}_0(D_j)$
be such that $\chi_j=1$ on $D_j\cap \phi^{-1}((-\infty, \phi(U_0)-\eps_0))$ for $\eps_0>0$ fixed but arbitrarily small.
If
$$
f_j=h^{-n/4} c_j(h) e^{-\frac{1}{h}\left(\phi(x)-\phi(U_j)\right)}\chi_j(x),\quad j=\pm 1,
$$
where $c_j(h)>0$ is a normalization constant such that $\norm{f_j}=1$, and
$$
\Pi=\frac{1}{2\pi i}\int_{\gamma} (z-P)^{-1}\,dz,\quad \gamma=\partial
D(0,\frac{h}{C}),\quad C>0,
$$
is the rank 2 spectral projection of $P$ corresponding to the eigenvalues
$\mu_0=0$ and $\mu_1$ in (\ref{eq3.1.1}), we have the basis
$$
e_j =\Pi f_j,\quad j=\pm 1,
$$
for ${\rm Ran}(\Pi)$, introduced in~\cite{HeHiSj}. From section 11
of~\cite{HeHiSj} we recall that
$$
e_j=f_j+{\cal O}(h^{-N_1}e^{-\frac{1}{h}(S_j-\eps_0)})\quad
\wrtext{in}\,\, L^2,\quad N_1>0,
$$
and that the restriction of $P$ to the space ${\rm Ran}(\Pi)$ has the matrix
\begeq
\label{eq3.2}
\left(
\begin{array}{cc}
\lambda_{-1}^*\\
\lambda_{1}^*
\end{array}
\right)
\left(
\begin{array}{cc}
\lambda_{-1} & \lambda_1
\end{array}
\right)=
\left(
\begin{array}{cc}
\lambda_{-1}^*\lambda_{-1} & \lambda_{-1}^* \lambda_1\\
\lambda_1^* \lambda_{-1} & \lambda_{1}^*\lambda_1
\end{array}
\right),
\endeq
with respect to the basis $(e_{-1}, e_1)$, with the eigenvalues $\mu_0=0$ and
$$
\mu_1=\lambda_{-1}^*\lambda_{-1}+\lambda_1^*\lambda_1.
$$
A simple computation shows that a corresponding basis of the
eigenvectors is given by
\begeq
\label{eq3.2.3}
\lambda_1 e_{-1}-\lambda_{-1} e_1
\endeq
and
\begeq
\label{eq3.2.4}
\lambda_{-1}^*e_{-1}+\lambda_1^* e_1.
\endeq
Here we recall from the formulas (11.43), (11.45), and the following discussion in~\cite{HeHiSj} that if
$\abs{\lambda_1}\geq \frac{1}{C}\abs{\lambda_{-1}}$ then $\abs{\lambda_1^*}\geq \frac{1}{2C}
\abs{\lambda_{-1}^*}$ and $\lambda_1\lambda_1^*>0$. We have the same fact after permuting the indices
$-1$, $1$ and the $\lambda_j$, $\lambda_j^*$. It follows that
\begeq
\label{eq3.2.5}
\mu_1\sim {\rm max} \abs{\lambda_j}^2\sim {\rm max}
\abs{\lambda_j^*}^2.
\endeq

\medskip
Rather than using (\ref{eq3.2.3}) and (\ref{eq3.2.4}), we shall make a normalized choice of the eigenfunctions,
given by
\begeq
\label{eq3.3}
v_0=\frac{1}{\sqrt{\mu_1}}\left(\lambda_1 e_{-1}-\lambda_{-1}e_1\right),
\endeq
and
\begeq
\label{eq3.4}
v_1=\frac{1}{\sqrt{\mu_1}}\left(\lambda_{-1}^* e_{-1}+\lambda_1^* e_1\right).
\endeq
The corresponding matrix of the coefficients is given by
\begeq
\label{eq3.5}
V=\frac{1}{\sqrt{\mu_1}}
\left(
\begin{array}{cc}
\lambda_1 & -\lambda_{-1}\\
\lambda_{-1}^* & \lambda_1^*
\end{array}
\right).
\endeq
We have ${\rm det}\, V=1$ and it follows from (\ref{eq3.2.5}) that
$V={\cal O}(1)$. Hence the inverse matrix $V^{-1}$ has the same properties, so that
$v_0$, $v_1$ is a well-behaved basis of eigenfunctions for $P$. If
$(e_{-1}^*, e_1^*)\in {\rm Ran}(\Pi^*)$ is the basis that is dual to
$(e_{-1}, e_1)$, then the corresponding basis of eigenfunctions of
$P^*$, dual to $(v_0, v_1)$ is given by the matrix $^t V^{-1}$, so that
\begeq
\label{eq3.6}
v_0^*=\frac{1}{\sqrt{\mu_1}}\left(\lambda_1^*e_{-1}^* - \lambda_{-1}^* e_1^*\right),
\endeq
and
\begeq
\label{eq3.7}
v_1^* = \frac{1}{\sqrt{\mu_1}}\left(\lambda_{-1}e_{-1}^*+\lambda_1 e_1^*\right).
\endeq

\medskip
We summarize the discussion above in the following proposition.
\begin{prop}
Let $v_j$ and $v_j^*$ be defined as in {\rm (\ref{eq3.3})},{\rm (\ref{eq3.4})}, {\rm (\ref{eq3.6})}, {\rm (\ref{eq3.7})}.
Then the spectral projections
$$
\Pi_j = (\cdot|v_j^*)v_j,\quad j=0,1
$$
associated to the eigenvalues $\mu_0=0$ and $\mu_1$ in {\rm (\ref{eq3.1.1})} are uniformly bounded as $h\rightarrow 0$.
\end{prop}

\bigskip
Combining (\ref{eq27}) together with Proposition 3.1, as well as with
Theorem 8.4 of~\cite{HeHiSj}, we get the result in Theorem 1.1 in the general case.

\begin{theo} \label{t32}
Let $P=-\Delta_A^{(0)}$ where we assume {\rm (\ref{eq3.0}-\ref{eq3.01})},
{\rm (\ref{morse})}, and {\rm (\ref{eq3.1.01})}. We also
assume that $P$ satisfies the dynamical hypotheses {\rm (\ref{eq1.91})}, {\rm (\ref{eq1.12.5})}, {\rm (\ref{eq1.12.54})},
so that the disc $D(0,h/C)$ for $C>0$ large enough, contains precisely {\rm 2}
eigenvalues of $P$, $\mu_0=0$ and $\mu_1$ given in {\rm (\ref{eq3.1.1})}. Let $\Pi_j$ be the spectral projection
associated with the eigenvalue $\mu_j$, $j=0,1$. Then we have
\begeq
\label{eq3.8}
\Pi_j={\cal O}(1),\quad h\rightarrow 0.
\endeq
We have furthermore, uniformly as $t\geq 0$ and $h\rightarrow 0$,
\begeq
\label{eq3.9}
e^{-tP/h}=\Pi_0+e^{-t\mu_1/h}\Pi_1+{\cal O}(1) e^{-t/C}, \quad C>0,\quad \wrtext{in}\quad {\cal L}(L^2,L^2).
\endeq
\end{theo}
Here we have also used that the eigenvalues of $P$ in $D(0,Bh)\backslash D(0,h/C)$ have real parts $\geq h/{\cal O}(1)$.

\section{Tunnel effect for a well and the sea}
\setcounter{equation}{0} In this section we shall show how to adapt the
analysis of section 11 of \cite{HeHiSj} and that of section 3 of the present work
to cover the case of a potential with a single well and a saddle
point, rather than a double well and a saddle point as before.  Some parts  of
this section are very close to the corresponding ones of section 8
in~\cite{HeHiSj}, and rather than repeating the arguments, we shall
often merely refer to the discussion there.

As in section 3, we shall consider the supersymmetric case. Assume
that we are given the constant matrices $A=B+C$ and a Morse function $\phi$,
that satisfy (\ref{eq3.0})--(\ref{eq3.01}). We then have the
corresponding Witten-Hodge Laplacian in degree 0, given by
(\ref{eq3.02}), with a principal symbol (\ref{eq3.1})
so that the assumptions (\ref{eq1.5})--(\ref{eq1.7}) hold. We refer to
the formula (11.3) of~\cite{HeHiSj} for the more general expression
for the Witten-Hodge Laplacian in degree $q\geq 0$, $P^{(q)}$.

\medskip
As before, we shall assume that $\phi$ has finitely many critical points
$x_1, ..., x_N \in \real^n$ and that $\vert \phi ' (x)\vert\ge 1/C$ when
$\vert x\vert\ge C$, with $C$ large. The assumptions (\ref{eq1.91})
is therefore satisfied with ${\cal C}=\{ \rho _j;\, j=1,...,N\} $ where
$\rho _j=(x_j,0)$. We also assume that the dynamical assumptions
(\ref{eq1.12.5}), (\ref{eq1.12.54}), (\ref{eq1.12.546}) hold.

As in section 11 of~\cite{HeHiSj}, an application of Theorem 8.3 of
\cite{HeHiSj} to $P^{(q)}$ shows that the eigenvalues $\mu_{j,k}$ there  are of the form
\begin{equation} \label{2w.0}
\mu _{j,k}={1\over i}\sum_{k=1}^n \sep{ \nu
_{j,k,\ell} + {1\over 2}}\lambda _\ell +\gamma _{j,k},
\end{equation}
where $\gamma_{j,k}$ is any eigenvalue of the subprincipal symbol $S_{P^{(q)}}$
at $(x_j,0)$. From the calculations in subsection 10.3 of
\cite{HeHiSj} we recall that the $\mu _{j,k}$ will be confined to a
sector $\{ 0\}\cup \{ |{\rm arg\, }z\vert <{\pi /2}-1/C\}$ around
$[0,+\infty)$. Recall that it is precisely when $x_j$ is of index $q$
(i.e. when the Hessian of $\phi $ at $x_j$ has precisely $q$
negative eigenvalues) that one of the $\mu _{j,k}$ is equal to
$0$.

\par We shall now introduce more specific conditions for the case that
we study here. Instead of assuming that we are in the double well
case, let us shall suppose that we have a single well and a sea, that is
 \begin{equation}\label{2w.0.5}
 \begin{split} &\phi
\hbox{ has precisely two critical points, one local} \\
&\hbox{ minimum } U_{ 1}, \hbox{ and a "saddle point" }U_0\hbox { of
index one.}
\end{split}
\end{equation} Notice that this implies that the Maxwellian
$e^{-\phi/h}$ is no longer an eigenfunction of $P^{(0)}$, since $\phi (x)$ does not go to $+\infty$ with $\abs{x}$.

\par Put $S_1=\phi (U_0)-\phi (U_1)$  so that $S_1>0$. The set
$\phi^{-1} (]-\infty ,\phi (U_0)[)$ has precisely two connected
components $D_j$, $j=\pm 1$, where $D_1$ is determined by the condition
that $U_1\in D_1$, while $D_{-1}$ is unbounded.

Under these assumptions we shall prove the following result:
\begin{theo} \label{wellsea} \label{2w1}
Let $P=-\Delta_A^{(0)}$ be as in {\rm (\ref{eq3.02})}, where we assume {\rm (\ref{2w.0.5})}.
Then for $C>0$ large enough,  $P$
has precisely {\rm 1} eigenvalue $\mu _1$  in the disc $D(0,h/C)$
when $h>0$ is small enough. Here $\mu _1$ is real and of the form
\begin{equation} \label{2w.28}  \mu _1=ha_1(h)e^{-2S_1/h},
\end{equation}
where $a_1(h)$ are real, $a_1(h)\sim a_{1,0}+a_{1,1}h+...$,
$a_{1,0}>0$, $S_1=\phi (U_0)-\phi (U_1)$.
\end{theo}

\remark It is clear that Theorem \ref{wellsea} implies an analog of Theorem 3.2 in the present metastable case. We
shall refrain from formulating it explicitly.

\begin{figure}[hbtp]
\begin{center}
\begin{picture}(0,0)%
\includegraphics{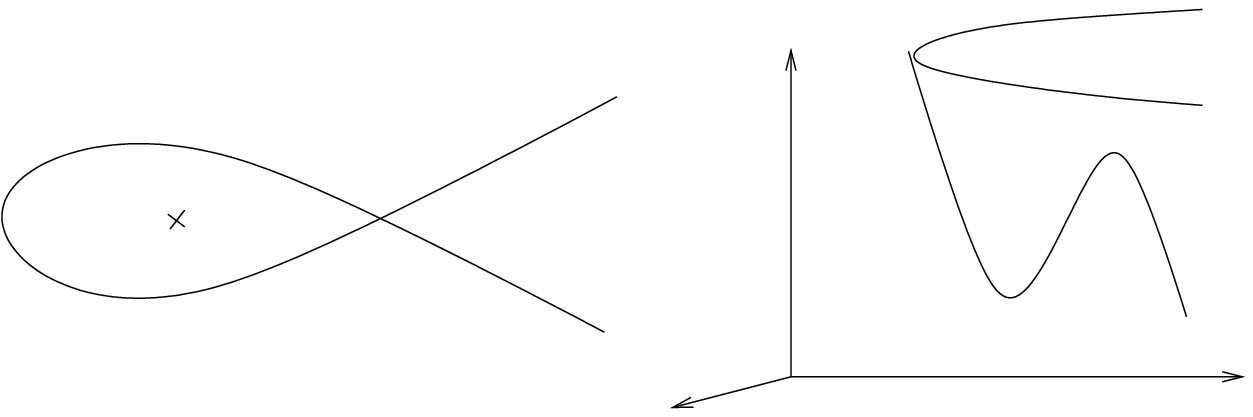}%
\end{picture}%
\setlength{\unitlength}{3158sp}%
\begingroup\makeatletter\ifx\SetFigFont\undefined%
\gdef\SetFigFont#1#2#3#4#5{%
  \reset@font\fontsize{#1}{#2pt}%
  \fontfamily{#3}\fontseries{#4}\fontshape{#5}%
  \selectfont}%
\fi\endgroup%
\begin{picture}(7483,2823)(655,-2215)
\put(1838,-849){\makebox(0,0)[lb]{\smash{{\SetFigFont{10}{12.0}{\rmdefault}{\mddefault}{\updefault}{$U_1$}%
}}}}
\put(1176,-937){\makebox(0,0)[lb]{\smash{{\SetFigFont{10}{12.0}{\rmdefault}{\mddefault}{\updefault}{$D_1$}%
}}}}
\put(3801,-837){\makebox(0,0)[lb]{\smash{{\SetFigFont{10}{12.0}{\rmdefault}{\mddefault}{\updefault}{$D_{-1}$}%
}}}}
\put(3351,-837){\makebox(0,0)[lb]{\smash{{\SetFigFont{10}{12.0}{\rmdefault}{\mddefault}{\updefault}{$U_0$}%
}}}}
\put(5106,454){\makebox(0,0)[lb]{\smash{{\SetFigFont{10}{12.0}{\rmdefault}{\mddefault}{\updefault}{$\phi(x)$}%
}}}}
\put(7264,-2149){\makebox(0,0)[lb]{\smash{{\SetFigFont{10}{12.0}{\rmdefault}{\mddefault}{\updefault}{$U_0$}%
}}}}
\put(6651,-1886){\makebox(0,0)[lb]{\smash{{\SetFigFont{10}{12.0}{\rmdefault}{\mddefault}{\updefault}{$U_1$}%
}}}}
\put(1364,-1699){\makebox(0,0)[lb]{\smash{{\SetFigFont{10}{12.0}{\rmdefault}{\mddefault}{\updefault}{level line at $\phi(U_0)$}%
}}}}
\end{picture}%
\caption{A well and the sea}
\end{center}
\end{figure}

\begin{proof} We first know that $P^{(0)}=-\Delta _A^{(0)}$ has
precisely one eigenvalue $\mu _1=o(h)$ spanning a corresponding
1-dimensional spectral subspace $E^{(0)}$ since there is a unique local minimum for $\phi$. Now
while $e^{-\phi/h}$ does not belong to $L^2$, a truncation of this function can be used as a quasimode near
$U_{ 1}$ and it follows therefore as in~\cite{HeHiSj}, that $\mu _{1}={\cal O}(h^\infty )$.
Moreover, $-\Delta _A^{(1)}$ has precisely one eigenvalue $\widetilde{\mu }_1=o(h)$ and $-\Delta _A^{(k)}$
has no eigenvalue $=o(h)$ for $k\ge 2$ from the discussion in the
beginning of the paragraph. Since our operators are real we know
that the spectra are symmetric around the real axis, hence $\mu
_1,\widetilde{\mu }_1$ are real. From the intertwining relations
$$
-\Delta _A^{(1)}d_\phi =d_\phi (-\Delta _A^{(0)}),\ -\Delta
_A^{(0)}d_\phi ^{A,*}=-d_\phi ^{A,*}\Delta _A^{(1)},
$$
we then also know that $\widetilde{\mu }_1=\mu _1$ (see also the
discussion at the end of page 69 of \cite{HeHiSj}).

\par The construction of the eigenfunctions $e_0$ and $e_1$ associated
to the critical points $U_1$ and $U_0$  for respectively
$-\Delta _A^{(1)}$ and $-\Delta _A^{(0)}$ is exactly the same as in \cite{HeHiSj}. We only retain the following from there:

\par We begin with $-\Delta _A^{(0)}$. Let $\chi _1\in C_0^\infty (D_1)$ be equal to 1 on
$D_1\cap\phi^{-1} (]-\infty ,\phi (U_0)-\epsilon _0])$ for
$\epsilon _0>0$ fixed but arbitrarily small. Consider
\begin{equation}\label{2w.1}
 f_1=h^{-n/4}c_1(h)e^{-{1\over h}(\phi
(x)-\phi (U_1))}\chi _1(x),
\end{equation}
 where $c_1(h)\sim c_{1,0}+hc_{1,1}+...>0$ is a normalization constant with
$c_{1,0}>0$, such that $  \Vert f_1\Vert =1.$
Then the normalized eigenfunction associated to $\mu_1$ is given by
\begin{equation} \label{2w.4}  e_1:=f_1+{\cal
O}(h^{-N_1}e^{-{1\over h}(S_1-\epsilon _0)})\hbox{ in }L^2.
\end{equation}

\par We continue with the study of $-\Delta _A^{(0)}$. Let
 $E^{(1)}$ be the one-dimensional eigenspace of $P^{(1)}$
corresponding to $\mu _1$. From an easy extension of
\cite[Theorem 9.1]{HeHiSj} (see also Remark 9.2 there) to the
non-scalar case with the presence of the other non-resonant well
$U_{1}$, we know that $E^{(1)}$ is generated by an eigenform
\begin{equation} \label{2w.10}  e_0(x;h)=\chi _0(x)e^{-{1\over h}\phi
_+(x)}h^{-{n\over 4}}a_0(x;h)+{\cal O}(e^{-S_0/h}),
\end{equation}
where $\chi _0\in C_0^\infty ({\rm neigh\,}(U_0))$ is equal to one
near $U_0$, $S_0>0$, and
$$
a_0(x;h)\sim \sum_0^\infty  a_{0,k}(x)h^k
$$
is a symbol as in Theorem 9.1 of~\cite{HeHiSj}, with $a_{0,0}(U_0)\ne
0$. Here the phase $\phi _+\in C^\infty ({\rm neigh\,}(U_0);[0,\infty))$
satisfies
$ \phi _+(x)\backsim \vert
x-U_0\vert ^2 $ and solves the eikonal equation
$ q(x,\phi _+'(x))=0$, with
 $q=p_2+p_1-p_0$.

\par From~\cite{HeHiSj}, let us recall that the phase function $\phi_+$ arises as the
generating function for the stable outgoing manifold through $(U_0,0)$
for the $H_q$-flow, $\Lambda_{\Phi_+}$, and recall also that
$\phi _+''(U_0)>0$ by Proposition \cite[Proposition 8.2]{HeHiSj}.
(Similarly we have a stable incoming manifold $\Lambda _{\phi
_-}$.) Let $k_\pm$ be the number of eigenvalues of the
linearization of ${{H_q}_\vert }_{\Lambda _\phi }$ at that point
with $\pm$ real part $>0$, so that $k_++k_-=n$. Let $K_+,
K_-\subset \Lambda _\phi $ be the corresponding stable outgoing
and incoming submanifolds of dimension $k_+$ and $k_-$
respectively. Then $K_+\subset \Lambda _{\phi _+}$, $K_-\subset
\Lambda _{\phi _-}$ and $\phi -\phi (U_0)-\phi _\pm$ vanishes to
the second order on $\pi _x(K_\pm )$. Since $\phi ''(U_0)$ has
signature $(n-1,1)$, we conclude that ${\rm dim\,}K_+=n-1$, ${\rm
dim\,}K_-=1$.  It is also clear that $\Lambda _\phi ,\Lambda
_{\phi _\pm}$ intersect cleanly along $K_{\pm}$, so we get
\begin{equation} \label{2w.10.3}
 \begin{split}
 \phi _+-(\phi -\phi (U_0))\backsim {\rm
dist\,}(x,\pi _x(K_+))^2,\\
 \phi -\phi (U_0)-\phi _-\backsim {\rm
dist\,}(x,\pi _x(K_-))^2.
\end{split}
\end{equation}

\par We now make some remarks about the adjoint operator $-\Delta _{\trans A}=(-\Delta _A)^{\trans A,*}$. As proved in \cite[Subsection 10.4]{HeHiSj}), we get the corresponding
 phases and submanifolds in this case, which satisfy
\begin{equation} \label{2w.10.4}
\begin{split}
  \phi
_+^*-(\phi -\phi (U_0))\backsim {\rm dist\,}(x,\pi _x(K_+^*))^2,\\
\phi -\phi (U_0)-\phi _-^*\backsim {\rm dist\,}(x,\pi
_x(K_-^*))^2.
\end{split}
\end{equation}
and for symmetry reason we recall that
\begin{equation} \label{2w.10.6}  \phi
_-^*=-\phi_+,\ \phi _+^*=-\phi_- ,
\end{equation} giving in
particular from (\ref{2w.10.3}), (\ref{2w.10.4}),
\begin{equation}
\label{2w.10.7}
\begin{split}  \phi -\phi (U_0) +\phi _+^*\backsim {\rm
dist\,}(x,\pi _x(K_-))^2,\\ \phi -\phi (U_0)+\phi _+\backsim {\rm
dist\,}(x,\pi _x(K_-^*))^2.
\end{split}
\end{equation}

\par Let  $\mu _1^*$ be the  eigenvalue of $P_*^{(0)}:=-\Delta
_{\trans A}^{(0)}$ that is  $o(h)$. As before, this is also the eigenvalue $o(h)$
 of $-\Delta
_{\trans A}^{(1)}$ and the corresponding eigenspaces $E_*^{(0)}$
 and $E_*^{(1)}$ are respectively spanned by  the eigenfunctions
\begin{equation}
\label{2w.10.8}
\begin{split} e_0^*(x;h)=\chi _0(x)e^{-{1\over h}\phi
_+^*(x)}h^{-{n\over 4}}a_0^*(x;h)+{\cal O}(e^{-S_0/h}) \\
\textrm{ and } \ \ \
e_1^*(x,h) = h^{-n/4}c_j(h)e^{-{1\over h}(\phi
(x)-\phi (U_1))}\chi _1(x)+ {\cal
O}(h^{-N_1}e^{-{1\over h}(S_1-\epsilon _0)}) \\
\end{split}
\end{equation}

\par Now, using that our eigenvalues and operators are real, we know by
duality that
$\mu _1^*=\mu _1$,
and that $(E_*^{(0)},E^{(0)})$ and $(E_*^{(1)},E^{(1)})$ are dual
pairs for the scalar products $(u\vert v)_{L^2}$ and $(u\vert
v)_A$ respectively. Following  \cite[Subsection 10.3]{HeHiSj} we know that $  (a_{0,0}^*(U_0)\vert
a_{0,0}(U_0))_A\ne 0$ and that $e_0^*$ can be normalized so that
\begin{equation} \label{2w.12}  (e_0^*\vert e_0)_A=1. \end{equation}
Similarly, denoting by $e_1^*$ the $L^2$ normalized eigenfuncion
spanning $E_*^{(0)}$, we have from (\ref{2w.4}-\ref{2w.10.8})
\begin{equation}
\label{2w.13}  (e_1^*\vert e_1)=1+{\cal O}(e^{-{1\over Ch}}).
\end{equation}

\par Let $(\lambda _1)$ be the (scalar) matrix of $d_\phi
:E^{(0)}\to E^{(1)}$ with respect to the bases $\sep{e_1}$ and
$(e_0)$.  Let also
$
\sep{ \lambda _1^*}
$
be the (scalar valued) matrix of $d_\phi ^{A,*}$ for the same
bases. The eigenvalue $\mu _1$ can be viewed as the scalar $d_\phi
^{A,*}d_\phi :E^{(0)}\to E^{(0)}$.  We get
\begin{equation}
\label{2w.14}  \mu _1=\lambda _{1}^*\lambda _{1},
\end{equation}
and
\begin{equation}
\label{2w.15}  \overline{\lambda }_1=(e_0^*\vert d_\phi e_1)_A,\ \
\ \ \
 \overline{\lambda }_1^*=(g_1\vert
d_\phi^{A,*} e_0)_A,\ j=\pm 1, \end{equation}
  where
$g_1=e_1^*(1+{\cal O}(e^{-{1\over Ch}}))$
 is the vector in $E_*^{(0)}$ that is dual to
$e_1$. Here the complex conjugate signs are superfluous since we
work with real operators, eigenvalues  and functions.

\par We skip the computation of $\lambda_1$, which is exactly the same as in \cite{HeHiSj},  just recalling that
the main term is equal to
\begin{equation} \label{2w.20}  -c_{1}(h)h^{1-{n\over 2}}\int
\chi _{1}(x)\langle A(x)a_0^*(x;h)\vert d\chi (x)\rangle
e^{-{1\over h}(\phi _+^*(x)+\phi (x)-\phi (U_{1}))}dx.
\end{equation}
and can be evaluated thanks to the stationary phase using (\ref{2w.10.7}). We get
\begin{equation} \label{2w.24}  \lambda _{1}=h^{1\over
2}\ell_{1}(h)e^{-{1\over h}S_1}(1+{\cal
O}(e^{-{1\over Ch}})), \ \ \  \ell_{1}\sim
\ell_{1,0}+h\ell_{1,1}+...,\ \ \  \ell_{1,0}\ne 0.
\end{equation}
similarly, $\lambda_1^*$
 ca be evaluated in a dual point of view as in \cite{HeHiSj} and we also get
\begin{equation} \label{2w.27}
 \lambda _1^* =  h^{1\over
2}\ell_{1}^*(h)e^{-{1\over h}S_{1}}(1+{\cal O}(e^{-{1\over Ch}})), \ \ \ \ell_1^*(h)\sim
\ell_{1,0}^*+h\ell_{1,1}^*+...,\ \ \  \ell_{1,0}^*\ne 0.
\end{equation}

We eventually claim that $\ell_{1,0}\ell_{1,0}^*>0$. Indeed, this
number is real and different form zero and if we deform our
matrices to reach the selfadjoint case (with $A>0$) we see that we
have a positive sign). Combining this with \Ref{2w.14}, the proof
of Theorem 4.1 is complete.
\end{proof}

\section{Some models of KFP type operators}
\setcounter{equation}{0}

\subsection{Probabilistic description} Here we shall give some examples of
Kramers-Fokker-Planck type operators. We begin with a very short review of stochastic
calculus in order to explain their probabilistic origin,
and refer to the books \cite{Ne67}, \cite{Oks00} for more details. Part of
this material can be also found in~\cite{EPR99}, \cite{EH00} and
\cite{EH02}, from where the example of the chain of anharmonic
oscillators is taken.

Let $x(t) \in \real^n$ be a stochastic process satisfying the
following stochastic differential equation
\begin{equation}
\label{SDE}
dx(t) = b(x(t))  dt + \sigma dw,
\end{equation}
where $w$ is the $m$-dimensional Wiener process, $\sigma$ is a
linear map from $\real^m$ to $\real^n$, and $b$ is a $C^{\infty}$-vector
field on $\real^n$, all of whose derivatives are bounded. Under these
assumptions, there exists a unique global solution $x(t)$ of (\ref{SDE}), for a given initial data $x(0)= x$,
independent of $w$, in an adapted stochastic $L^2$ setting --- see the
references already mentioned. Then we can define a semigroup of
operators $T^t$, $t\geq 0$, by
\begin{equation}  \label{semi}
\E \sep{ \varphi(x(t)) | \fff^s } = T^{t-s} \varphi(x(s)), \ \ \ a.s.
\end{equation}
when $0 \leq s \leq t$.
Here $\fff^t$ is the filtration associated to $\set{ w(s) - w(0); \ 0
\leq s \leq t}$ and $x$, and $\varphi \in C_{(0)}(\real^n)$, where
$C_{(0)}(\real^n)$ is the Banach space of continuous functions vanishing at infinity, with the
topology of the uniform convergence. Then $T^t$ is a strongly continuous
positivity preserving contraction semigroup, whose infinitesimal
generator is given on $C^\infty_0(\real^n)$ by
$$
L = \nabla \cdot D \nabla + b(x)\cdot \nabla,
$$
where $D = \frac{ 1}{2} \sigma \sigma^t$.  The idea now is to extend
$T^t$ to a larger class of test functions, and then to study the evolution
of the adjoint $(T^t)^*$ on the dual space of bounded measures. To be precise,
let us denote by  $d\mu_t(x)$ the probability distribution for $x(t)$, defined for
all $t \geq 0$. We then have
$$
\E \sep{ \varphi(x(t)) } = \int \varphi(x) d\mu_t(x), \ \ \ \ \varphi
\in C_{(0)}(\real^n),
$$
and we get by (\ref{semi}) that $\mu_t = (T^t)^* \mu_0$, where
$(T^t)^*$ is the adjoint of the operator $T^t$ acting on the Banach space of
bounded measures on $\real^n$.

We shall now extend the space of test functions. When doing so, we introduce the Hilbert space
\begin{equation} \label{hhh}
\hhh = L^2(\real^n, e^{-\Phi(x)} dx)
\end{equation}
where $\Phi\in C^{\infty}(\real^n)$. We shall make the following
assumptions concerning $\Phi$ :
\begeq
\label{eq3.01bis}
\partial_x^{\alpha} \Phi(x)={\cal O}(1),\quad \abs{\alpha}\geq 2,
\endeq
\begeq
\label{eq3.01bis2}
\frac{1}{2} b(x)\cdot \nabla\Phi(x) +
\frac{1}{4} \langle{D\nabla_x\Phi,\nabla_x\Phi\rangle}-\frac{1}{2}{\rm div}\, b \leq {\cal O}(1)
\endeq
and
\begeq
\label{eq3.0100}
\partial_x^{\alpha} \bigg( \frac{1}{2} b(x)\cdot \nabla\Phi(x) +
\frac{1}{4} \langle{D\nabla_x\Phi,\nabla_x\Phi\rangle}-\frac{1}{2}{\rm div}\, b \bigg)={\cal O}(1),\quad \abs{\alpha}\geq 2.
\endeq
These conditions will be fulfilled in the case that we shall study in what
follows, since in the supersymmetric case it is straightforward to verify that they are equivalent to (\ref{eq3.01}). Now we can identify the dual $\hhh'$ of $\hhh$ with the space of densities
$$
\hhh^* = L^2(\real^n, e^{\Phi(x)} dx).
$$
Assume that the measures $d\mu_t$ are absolutely continuous with density in $\hhh^*$, and write
$$
d \mu_t = f(t,.) dx,
$$
identifying the measure $d\mu_t$ with the corresponding density $f_t$. We denote again
by $(T^t)^*$ acting on $\hhh^*$ the adjoint of $T^t$ acting on $\hhh$. We introduce the
formal adjoint operator $L^*$ on $\hhh^*$ of $L$ on $\hhh$, with the domain
$C^\infty_0(\real^n)$, which is given by
$$
L^* = \nabla\cdot D \nabla - \nabla \cdot b(x).
$$
We have the following result.

\begin{lemma}
Assume that $\Phi\in C^{\infty}(\real^n)$ satisfies {\rm
(\ref{eq3.01bis})}, {\rm (\ref{eq3.0100})}. Then  operator $-L$ on $\hhh$ (resp. $-L^*$
 on $\hhh^*$) is  m-accretive, and
$T^t$ (resp $(T^t)^*$) is a  strongly continuous semigroups on $\hhh$ (resp. $\hhh^*$), with infinitesimal generators given by $L$ (resp. $L^*$).
\end{lemma}
\preuve
It will be more convenient to work in the unweighted space $L^2(\real^n)$. To this
end, if $\phi\in {\cal H}$, we write $\phi(x) = e^{\Phi(x)/2} \psi(x)$, $\psi\in L^2$. If
$$
\partial_t \phi = L\phi,
$$
then the equation satisfied by $\psi$ is
$$
\partial_t \psi = \left(e^{-\Phi/2}L e^{\Phi/2}\right)\psi,
$$
so that
$$
\D_t \psi =  (\D_x + \D_x\Phi/2)\cdot D (\D_x + \D_x\Phi/2) \psi +  b(x)\cdot (\D_x + \D_x\Phi/2) \psi.
$$
Let $C$ be a sufficiently large constant. According to (\ref{eq3.01bis}-\ref{eq3.0100}),
the operator
$$
- \sep{ (\D_x + \D_x\Phi/2) \cdot D (\D_x + \D_x\Phi/2)  +  b(x)\cdot
(\D_x + \D_x\Phi/2)} + C
$$
has a symbol satisfying the hypotheses (\ref{eq1.4})--(\ref{eq1.7}). Here we may recall that the vector field $b$ is
bounded on $\real^n$ together with all of its derivatives. An application of Corollary 3.2 in~\cite{HeHiSj} shows that
its maximal closed realization in $L^2$ coincides with the graph closure on
${\cal S}$. Coming back to $\hhh$ and denoting by $L$ again its maximal closed extension,
we get that $L+C$ is maximal accretive, and that $T^t$ is a strongly
continuous semigroup, thanks to the Hille-Yosida Theorem. As for  the
dual semi-group on $\hhh^*$, we also get (for example, using Corollary
10.6 in \cite{Paz83}), that the same occurs for $(T^t)^*$ and
     $L^*$. The proof is complete.
\fin

From the preceding discussion, we get the equation
satisfied by the density $f$ for an initial data $f_0 \in \hhh^*$,
\begin{equation} \label{eqf}
\left\{ \begin{array}{l} \D_t f + (- L^*) f = 0 \\ f|_{t=0} = f_0
\end{array} \right. \ \ \ \ \mbox{ i.e.}  \ \ \
\left\{ \begin{array}{l} \D_t f + ( -\nabla \cdot D \nabla  +
\nabla \cdot b) f = 0 \\
f|_{t=0} = f_0
\end{array} \right.
\end{equation}
where we recall that $D = \sigma \sigma^t/2$. In particular we have
 $d \mu_t = f(t,.) dx = (T^t)^* \mu_0$ in the space of  bounded measures.

\par If there exists an
invariant probability measure  $\mu_\infty$,
then its density $\mmm$ is in $\hhh^*$. In our present
study we shall essentially make the choice $C^{-1}e^{-\Phi} = \mmm$,
but there are cases (see e.g. \cite{EPR99}, \cite{EH00}), where it may happen that no
invariant measure is known, and that another choice of the function
$\Phi$ is necessary.  Such a function $\mmm$ will be called a Maxwellian
of the process. Notice that if it  exists, it is a $0$-eigenfunction
of $L^*$ and positive.

\remark Notice that it may also happen that there exists an invariant
measure, which fails to be finite. We also associated to it a function that we will call again
Maxwellian (and denote again by $\mmm$). In that case of course it
cannot be  normalized.
\bigskip

Equation (\ref{eqf}) is nearly the Kramers-Fokker-Planck type equation
that we studied in the first part of the paper.
In the following sections we shall also do the following two things:
first we shall exhibit the semiclassical scaling, which corresponds to
the low temperature limit in the models we are going to study later. Second,
we shall change our unknown by posing $f = e^{-\Phi(x)/2} u$ (forgetting for a while the semiclassical scaling),  in order
to work in the flat space $L^2$ rather than in $\hhh$. Finally,
in the three models that we present in the next subsections (Witten, Kramers-Fokker-Planck, and
the chain of anharmonic oscillators), we shall  recognize the
supersymmetric structure.

\subsection{Witten and Kramers-Fokker-Planck operators} We begin with the
Witten case. It corresponds to an evolution equation with a gradient
field $ -\gamma \nabla V(x)$ and a diffusion force coming from a heat
bath at a temperature $T$. We have
$$
dx = - \gamma \D_x V dt + \sqrt{2\gamma T} d w .
$$
Here $x\in \real^n$ is the spatial variable, the parameter $\gamma$ is a friction coefficient, and $w$ is an
$n$--dimensional Wiener process of mean $0$ and variance $1$.  With
the notation of the preceding subsection, we recover an equation of type (\ref{SDE}) with
$D = \sigma \sigma^t /2 = \gamma T I_d$ and $b(x) = -\gamma \D_x
V$. Equation (\ref{eqf}) for the density in this case is then
\begin{equation} \label{bef1}
\left\{ \begin{array}{l} \D_t f  - \gamma T \D_x^2 f - \gamma \D_x
(\D_x V f)  = 0 \\ f|_{t=0} = f_0
\end{array} \right. \ \ \ \ \mbox{ i.e.}  \ \ \
\left\{ \begin{array}{l} \D_t f -\gamma \D_x ( T \D_x + \D_x V) f = 0 \\
f|_{t=0} = f_0.
\end{array} \right.
\end{equation}
Posing $T=h/2$ and multiplying by $h$ gives the semiclassical equation
\begin{equation}
\label{W1}
 h\D_t f -\frac{\gamma}{2} h\D_x (h\D_x + 2\D_x V) f = 0.
 \end{equation}
It is then clear that an associated Maxwellian of the process is
$$
\mmm(x) = e^{-2V(x)/h}.
$$
Writing $f = \mmm^{1/2} u$, we obtain from (\ref{W1}) that
\begin{equation} 
 h\D_t u + \frac{\gamma}{2} ( - h\D_x + \D_x V)(h\D_x + \D_x V) u = 0.
 \end{equation}
Here we recognize the Witten operator $W = ( - h\D_x + \D_x V)
(h\D_x + \D_x V)$. In the notation of (\ref{eq3.0})--(\ref{eq3.02}),
it corresponds to a supersymmetric operator with
$$
A= \frac{\gamma}{2} I_d, \ \ \ \phi(x) = V(x)
$$
Assumptions of type (\ref{eq3.01}) on $V$ are then fulfilled if
$$
\D^\alpha V(x)  = \left\{\begin{array}{l} \ooo(1) \ \ \textrm{ when
} \ \ |\alpha| =2
\\ \ooo( \langle{x\rangle}^{-1}) \ \  \textrm{ when } |\alpha| \geq 3
\end{array} \right.
$$
If we also suppose that $V$ is a Morse function with two local minima and a
saddle point of index one, such that
\begeq
\abs{\nabla V}\geq 1/C,\quad \wrtext{for}\quad \abs{x}\geq C>0.
\endeq
then the dynamical asumptions (\ref{eq1.12.5})--(\ref{eq1.12.546}) are
 satisfied (we skip the proof here, which will be given later in the more complex case of the chain of oscillators). In particular, in this case we get Theorem \ref{t32}.
Of course, the corresponding result follows also in the case of a
single well and the sea, i.e. when $V$ has precisely one local minimum and a saddle point (Theorem
\ref{2w1} and the remark following it).

\bigskip
\par
We proceed now to discuss the Kramers-Fokker-Planck case, and follow  the
same method. The stochastic equation of type (\ref{SDE}) comes here from the
Newton law
\begin{equation}
\left\{ \begin{array}{l} dx = y dt \\ dy =  - \gamma y dt- \D_x V(x) dt
+ \sqrt{2\gamma T} d w \end{array} \right..
\end{equation}
The parameter $\gamma$ is a friction coefficient, and the particle of
position $x \in \real^n$ and velocity $y\in \real^n$ is submitted to an external force
field derived from a potential $V$, with $w$ being an $n$--dimensional
Brownian process of mean $0$ and variance $1$. With the notation of
the preceding subsection, we therefore have
 $$
 D =  \sigma \sigma^t /2 = \left( \begin{array}{cc} 0 & 0 \\ 0 & \gamma T I_d \end{array}\right)
\ \ \ \ \textrm{ and } \ \ \ \  b(x,y) = \left( \begin{array}{cc} y \\ - \gamma
y - \D_x V
  \end{array}\right).
  $$
The corresponding equation for the density (\ref{eqf}) is then
\begin{equation} 
\left\{ \begin{array}{l} \D_t f  - \gamma T \Delta_y f + \D_x (y f)
+ \D_y ( -\gamma y  f -  \D_x V f )  = 0 \\ f|_{t=0} = f_0
\end{array} \right.
\end{equation}
\mbox{ i.e.}  \ \ \
\begin{equation} \label{bef2}
\left\{ \begin{array}{l} \D_t f -\gamma \D_y .  ( T \D_y + y) f +
 y \D_x f -  \D_x V  \D_y  f  = 0 \\
f|_{t=0} = f_0.
\end{array} \right.
\end{equation}
Posing $T=h/2$ and multiplying by $h$ gives the semiclassical formulation
\begin{equation}
\label{KFP}
 h\D_t f -\frac{\gamma}{2} h\D_y .  ( h \D_y + 2y) f +
 y h\D_x f -  \D_x V  h\D_y  f  = 0.
 \end{equation}
A Maxwellian of the process is then
$$
\mmm(x,y) = C^{-1} e^{-2(V(x) + {y^2}/2)/h}
$$
where $C$ is a normalization constant. If we write $f = \mmm^{1/2} u$,
then (\ref{KFP}) gives
\begin{equation} 
h\D_t u +\frac{\gamma}{2}(-h\D_y + y)\cdot ( h \D_y + y) u +
\gamma y \cdot h\D_x u -  \D_x V \cdot h\D_y  u  = 0
\end{equation}
This is the
semiclassical Kramers-Fokker-Planck operator (\ref{eq01}) that we
studied in \cite{HeSjSt}, \cite{HeHiSj}, and the first part of the
present paper. In the notation of section 3, it corresponds to a supersymmetric operator with
$$
 A =  \frac{1}{2} \left( \begin{array}{cc} 0 & I_d \\ - I_d & \gamma  \end{array}\right)
\ \ \ \ \textrm{ and } \ \ \ \  \phi(x,y) = V(x) +
y^2/2.
$$
The  assumptions (\ref{eq3.01}) are fulfilled if
$$
\D^\alpha V(x)  =  \ooo(1) \ \ \textrm{ when
} \ \ |\alpha|  \geq 2
$$
As in the Witten case, if we also suppose that $V$ is a Morse function
with precisely two local minima and a saddle point of index one, and that
\begeq
\abs{\nabla V}\geq 1/C,\quad \wrtext{for}\quad \abs{x}\geq C>0.
\endeq
then the dynamical asumptions (\ref{eq1.12.5})--(\ref{eq1.12.546}) are
satisfied, and Theorem \ref{t32} is applicable. The corresponding
result in the case of a single well and the sea, i.e. when $V$ has one
local minimum and a saddle point, is also valid (Theorem \ref{2w1} and
the following remark).

\section{Chains of anharmonic oscillators}\label{section0}
\setcounter{equation}{0}

The last example that we give comes from the series of papers
\cite{EPR99}, \cite{EH00} \cite{EH02}. It is a model describing   a
chain of two anharmonic oscillators coupled with two heat baths at
each side.




The particles are described by their respective position and velocity
$(x_j, y_j) \in \real^{2d}$.
We suppose that for each oscillator $j \in \set{ 1,2}$, the particles
are  submitted to an external force derived from a potential
$V_j(x_j)$, and that there is a coupling between the two oscillators
derived from a potential $V_c(x_2-x_1)$.  We denote by $V$ the sum
$$
V(x) = V_1(x_1)+ V_2(x_2) + V_c(x_2-x_1),
$$
where $x = (x_1, x_2)$, and we also write $y = (y_1,y_2)$.
By $z_j$, $j \in \set{ 1,2}$ we shall denote the variables describing
the state of the particles in each of the heat baths, and set $z = (z_1,z_2)$. We suppose that the
particles in each bath are submitted to a coupling with the nearest
oscillator, a friction force  and a thermal diffusion at temperature
$T_j$, $ (j=1,2)$. We denote by $w_j$, $j \in \set{ 1,2}$, two
$d$-dimensional brownian motions of mean 0 and variance $1$, and set $w=(w_1, w_2)$.

\begin{figure}[hbtp]
\begin{center}
\begin{picture}(0,0)%
\includegraphics{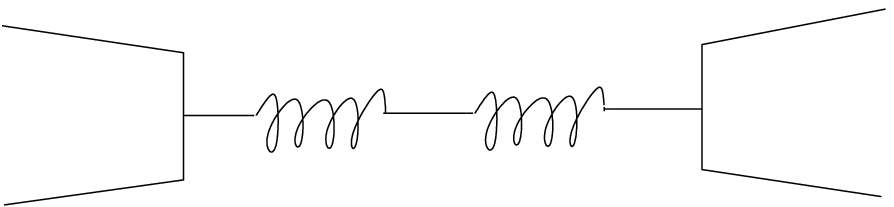}%
\end{picture}%
\setlength{\unitlength}{3158sp}%
\begingroup\makeatletter\ifx\SetFigFont\undefined%
\gdef\SetFigFont#1#2#3#4#5{%
  \reset@font\fontsize{#1}{#2pt}%
  \fontfamily{#3}\fontseries{#4}\fontshape{#5}%
  \selectfont}%
\fi\endgroup%
\begin{picture}(5324,1199)(164,-2348)
\put(1701,-1474){\makebox(0,0)[lb]{\smash{{\SetFigFont{10}{12.0}{\rmdefault}{\mddefault}{\updefault}{$(x_1,y_1)$}%
}}}}
\put(3051,-1461){\makebox(0,0)[lb]{\smash{{\SetFigFont{10}{12.0}{\rmdefault}{\mddefault}{\updefault}{$(x_2,y_2)$}%
}}}}
\put(726,-1886){\makebox(0,0)[lb]{\smash{{\SetFigFont{10}{12.0}{\rmdefault}{\mddefault}{\updefault}{$z_1$}%
}}}}
\put(339,-1636){\makebox(0,0)[lb]{\smash{{\SetFigFont{10}{12.0}{\rmdefault}{\mddefault}{\updefault}{$\alpha_1 h$}%
}}}}
\put(5026,-1899){\makebox(0,0)[lb]{\smash{{\SetFigFont{10}{12.0}{\rmdefault}{\mddefault}{\updefault}{$z_2$}%
}}}}
\put(4639,-1636){\makebox(0,0)[lb]{\smash{{\SetFigFont{10}{12.0}{\rmdefault}{\mddefault}{\updefault}{$\alpha_2 h$}%
}}}}
\end{picture}%
\caption{Oscillators coupled to heat baths}
\end{center}
\end{figure}

The fundamental system of equations of type (\ref{SDE}) is then written
as follows, (see
\cite{EPR99} for more detail concerning the physical constants)
\begin{equation}
\left\{
\begin{array}{l}
 dx_1 = y_1 dt \\
  dy_1 =  - \D_{x_1} V(x) dt + z_1 dt \\
  d z_1 = - \gamma z_1 dt + \gamma x_1 dt - \sqrt{2\gamma T_1} d
  w_1 \\
  d z_2 = - \gamma z_1 dt + \gamma x_2 dt - \sqrt{2\gamma T_2} d
  w_2 \\
dy_2 =  - \D_{x_2} V(x) dt + z_2 dt \\
d x_2 =y_2 dt.
\end{array} \right.
\end{equation}
The parameter $\gamma$  is  the  friction coefficient in the baths.
In the $(x,y,z)$ variables, the diffusion matrix and the drift
appearing in (\ref{SDE}) are
therefore
$$
 D =  \sigma \sigma^t /2 = \left( \begin{array}{ccc} 0 & 0 & 0 \\ 0&0&0 \\
 0 & 0 &  \gamma T I_d \end{array}\right)
\ \ \ \ \textrm{ and } \ \ \ \  b(x) = \left( \begin{array}{cc} y \\  - \D_x
V + z \\ \gamma(x-z)
  \end{array}\right).
  $$
  (for simplicity we identified  $T$  with the 2d times 2d diagonal matrix with
  coefficients $T_1 I$ and $T_2 I$.
The corresponding equation (\ref{eqf}) for the density  is then
\begin{equation} \label{bef5}
\left\{ \begin{array}{l} \D_t f  -  \gamma T \D_{z}^2 f  + \D_x (y
f)
+ \D_y ( -  \D_x V f + z f ) +  \D_z ( \gamma(x -z)) = 0 \\
f_{t=0} = f_0
\end{array} \right.
\end{equation}
where $T \D_z^2$ stands for $T_1 \D_{z_1}^2 + T_2 \D_{z_2}^2$. We
get
\begin{equation} \label{bef6}
\left\{ \begin{array}{l} \D_t f -\gamma \D_z  ( T \D_z + (z-x)) f
+  y \D_x f -  (\D_x V -z)  \D_y  f  = 0 \\
f_{t=0} = f_0.
\end{array} \right.
\end{equation}
Notice that it is very close to the Witten and Kramers-Fokker-Planck operator.
For a semiclassical formulation, we pose $T_1 = \alpha_1 h/2$ and $T_2 =
\alpha_2 h/2$, and we multiply (\ref{bef6}) by $h$. This gives
\begin{equation} \label{firstW6} \begin{split}
 h\D_t f & +  \frac{\gamma}{2} \alpha_1 (-h\D_{z_1})( h \D_{z_1} + 2(z_1-x_1)/\alpha_1)f \\
 & + \frac{\gamma}{2} \alpha_2 (-h\D_{z_2})( h \D_{z_2} + 2(z_2-x_2)/\alpha_2)f \\
  &  + (y h\D_x f -  (\D_x V -z)   h\D_y) f  = 0.
   \end{split}
 \end{equation}
 At this stage it is difficult to exhibit a Maxwellian. Indeed the
 existence of an invariant measure is a difficult problem solved in
 some particular case in \cite{EPR99}. Anyway, it is clear that the
 function
 $$
 \Phi(x,y,z) = V(x) + y^2/2 + z^2/2 - zx
 $$
 plays a special role, in fact it is the classical energy at
 temperature $1$. We can also check that in the case of same
 temperatures ($\alpha_1=\alpha_2 \defegal \alpha)$,
 the function
$$
\mmm_\alpha = C^{-1} e^{-2\Phi/\alpha h}
$$
is a Maxwellian of the process. We use this function to define
the weighted space $\hhh \defegal L^2(e^{-2\Phi/\alpha h} dxdydz)$ as in (\ref{hhh}), and in order to work in the flat space $L^2$ we make the change of  unknown
$$
f = \mmm_\alpha^{1/2} u.
$$
Equation (\ref{bef2}) becomes
\begin{equation} \label{firstW7} \begin{split}
 h\D_t u & +  \frac{\gamma}{2} \alpha_1 \sep{-h\D_{z_1} +
 \frac{1}{\alpha}(z_1-x_1)}
             \sep{ h \D_{z_1} + \sep{\frac{2}{\alpha_1}- \frac{1}{\alpha}}(z_1-x_1)}u \\
 & +  \frac{\gamma}{2} \alpha_2 \sep{-h\D_{z_2} + \frac{1}{\alpha}(z_2-x_2)}
             \sep{ h \D_{z_2} + \sep{\frac{2}{\alpha_2}- \frac{1}{\alpha}}(z_2-x_2)}u \\
  &  + \sep{y h \D_x  -  (\D_x V -z)   h \D_y} u  = 0.
   \end{split}
 \end{equation}
 We impose the following condition on the parameter $\alpha$ :
$$
 \alpha \geq \max \set{ \alpha_1, \alpha_2}/2.
$$
which  corresponds to a semiclassical study at "reference" temperature
$\alpha h/2$ not too low.

Unfortunately we are  not able to find any supersymmetric structure
 in the case of different temperatures, since  a Maxwellian
  is not known in this case. From now on we therefore stick to the
 case of identical temperatures $T=h/2$ so that
  $$
  \alpha = \alpha_1 = \alpha_2=1.
  $$
  Equation (\ref{firstW7}) becomes
  \begin{equation} \label{firstW8} \begin{split}
 h\D_t u & +  \frac{\gamma}{2}  \sep{-h\D_{z_1} +
 (z_1-x_1)}
             \sep{ h \D_{z_1} + (z_1-x_1)}u \\
 & +  \frac{\gamma}{2}  \sep{-h\D_{z_2} + (z_2-x_2)}
             \sep{ h \D_{z_2} + (z_2-x_2)}u \\
  &  + \sep{y h \D_x  -  (\D_x V -z)   h \D_y} u  = 0
   \end{split}
 \end{equation}
  and the Maxwellian  was already exhibited  $ \mmm_1 = C^{-1} e^{-2\Phi/ h}. $ This equation can be written $h \D_y u + P u = 0$ where
  \begeq \label{ana}
  \begin{split}
  P =&  \frac{\gamma}{2}  \sep{-h\D_{z_1} +
 (z_1-x_1)}
             \sep{ h \D_{z_1} + (z_1-x_1)} \\
 & +  \frac{\gamma}{2}  \sep{-h\D_{z_2} + (z_2-x_2)}
             \sep{ h \D_{z_2} + (z_2-x_2)}   + \sep{y h \D_x  -  (\D_x V -z)   h \D_y}.
  \end{split}
  \endeq
In the notations of Section 3 (\ref{eq3.0}-\ref{eq3.02}), we can write $P$ as a Witten-Hodge laplacian $P= -\Delta_A^{(0)}$ with a supersymmetric phase $\phi$  given by
\begin{equation} \label{phib}
\phi \defegal \Phi =  V(x) + y^2/2 + z^2/2 - zx,
\end{equation}
and the  non-degenerate matrix $A=B+C$   given by
$$
A = \frac{1}{2} \left( \begin{array}{ccc}
0 & I_d & 0 \\
-I_d & 0 & 0 \\
0 & 0 & \gamma I_d
\end{array}\right) \ \ \ \ \textrm{ with } \
B = \frac{1}{2} \left( \begin{array}{ccc}
0 &0 & 0 \\
0 & 0 & 0 \\
0 & 0 & \gamma I_d
\end{array}\right),  \ \ \ C = \frac{1}{2} \left( \begin{array}{ccc}
0 & I_d & 0 \\
-I_d & 0 & 0 \\
0 & 0 & 0
\end{array}\right).
$$
In order to
 complete the semiclassical study as in \cite{HeHiSj}, we only need
additional conditions on the potentials $V_1$, $V_2$ and $V_c$. It is clear
that the conditions
\begin{equation} \label{quad}
\D^\alpha V_\eps(x)  =  \ooo(1) \ \ \textrm{
when } \ \ |\alpha|  \geq 2, \ \ \ \ \textrm{ with }
 \eps = 1,2 \textrm{ and }c
\end{equation}
imply (\ref{eq3.01}).
  In view of
the definition (\ref{phib}), it is straightforward that $\phi$
has exactly the same number of  critical points than  $V(x)-x^2/2$ with same index.
 For this it is sufficient to notice that
there is a natural splitting of the variables for $\Phi$  given by
$$
 \Phi = (V(x) - x^2/2) + y^2/2 + (z-x)^2/2.
$$

 We postpone to the end of this section the proof of the following lemma:

\begin{lemma} \label{assume}
Suppose that $V$ satisfies {\rm (\ref{quad})}. If  in addition  $V(x) -x^2/2$
is a  Morse function and   there exists $C$ such that
\begin{equation} \label{tuc}
|\D V(x) -x| \geq 1/C    \textrm{ when } |x| \geq C,
\end{equation}
then {\rm (\ref{morse})} and  the dynamical conditions {\rm (\ref{eq1.12.5}-\ref{eq1.12.546})} are fulfilled.
\end{lemma}

As a consequence we can apply Theorem \ref{t32} to operator $P$:

\begin{prop}
Consider $P$ given by {\rm (\ref{ana})} and suppose that $V$ satisfies {\rm (\ref{quad})} and {\rm (\ref{tuc})}. Then if the
effective potential $V(x) - x^2/2$ is of double well type,
 (two local minima and a saddle point of index 1),  the hypotheses of Theorem {\rm \ref{t32}} are fulfilled and as a consequence its conclusions apply to operator $P$.
 \end{prop}

 \preuve It is straightforward.
 From the construction of $P$, hypotheses (\ref{eq3.0}-\ref{eq3.002}) are fulfilled. From (\ref{quad}), hypothesis  (\ref{eq3.01}) is satisfied. Since the
 effective potential $V(x) - x^2/2$ is a morse function of double well type, then  (\ref{eq3.1.01}) is also satisfied since, as
 already noticed,
 $V(x) - x^2/2$ and $\phi$ have the same number of critical point with same index. Eventually using Lemma
 \ref{assume} and (\ref{tuc}) we get that hypothese (\ref{morse}) and
 the dynamical conditions {\rm (\ref{eq1.12.5}-\ref{eq1.12.546})} are fulfilled. The proof is complete and  the  conclusions of Theorem \ref{t32} apply to $P$.
 \fin

 \remark
Of course the corresponding result  follows in the case of a well and the sea, ie when
$V(x) - x^2/2$ has one minimum and a saddle point  (Theorem \ref{2w1} and the remark after).

\bigskip
A simple family of such potentials is given for example by the ones for which $V_1(x_1)-x_1^2/2$ of double well
type, $V_2(x_2)-x_2^2/2$ of single well type, and $V_c$ sufficiently small.
Here is an example of such potentials in 1d:
$$
V_1(x_1) = x_1^2/2+ 5\sqrt{(x_1^2 - 1)^2 + 1}, \ \ \ V_2(x_2) = 5 x_2^2, \ \ \
V_c(x') = \frac{1}{10} \cos(x').
$$
Here $x_1, x_2, x' \in \real$.

\remark We did all the computations in the case of 2 oscillators. It
is clear that  the preceding supersymmetric construction works as well in the case of
$N$ oscillators coupled and with two heat bathes at each side as in
\cite{EPR99}.  We did not try to verify the dynamical conditions in these cases.
Eventually recall  that the complete study at  different temperatures seems difficult to treat
 (see e.g. the recent work by  Hairer and Mattingly \cite{HM07} in the case of 3 oscillators).

\preuve[of Lemma \ref{assume}]

We only prove here condition (\ref{morse}) and (\ref{eq1.12.546}) since the proof of  other ones follow the same kind of arguments.
The proof uses strongly the fact that the symbols are with quadratic growth at most.
We write  $p = p_2 + i p_1 + p_0$ for the symbol of the corresponding Hodge Laplacian on $0$-forms (minus the constant $\gamma h d/2$) where
denoting $(\xi, \eta, \zeta)$ the dual variables of $(x,y,z)$ we have
$$
p_2 = \frac{\gamma}{2} \zeta^2, \ \ \  p_1 =  y\xi - (\D_xV -z)\eta, \ \ \
p_0 = \frac{\gamma}{2} (z-x)^2.
$$
In particular,   with the notations of Section 2, we have
\begin{equation} \label{nu}
\nu(x,y,z, \D_x, \D_y, \D_z) =  y \D_x  - (\D_xV -z)\D_y.
\end{equation}
 We denote by $x_1, ..., x_N$ the critical points of $V(x) - x^2/2$, and notice that the critical points of $\Phi$ are  $(x_j, 0, x_j) $ for ${j=1, ...,N}$. According to definitions
 (\ref{eq1.10}) and (\ref{nu}), the critical set $\ccc$ of $p$ is made of the points $\rho_j = (x_j, 0, x_j, 0, 0, 0)$ for ${j=1, ..., N}$.
We also introduce $\pi_x$ (resp. $\pi_{xyz}$) the orthogonal  projections on $\R^d_x$ (resp. $\R^{3d}_{x,y,z}$) from $\R^{6d}$.

\bigskip Let now $\eps >0$ be a fixed constant.
 Since $V(x) - x^2/2$ is a Morse function, it is non-degenerate, so that with (\ref{tuc}) we get  the following:  there exists $C_\eps$ such that
\begin{equation} \label{tuc2}
\forall x \not\in \pi_x \sep{\ccc + B(0,\eps)}, \ \ \ \
|\D V(x) -x| \geq 1/C_\eps.
\end{equation}
From this result we get easily  the following one: there exists $C'_\eps$ such that
\begin{equation} \label{tuc3}
\forall (x,y,z) \not\in \pi_{x,y,z} \sep{\ccc + B(0,\eps)}, \ \ \ \
|\D V(x) -x| + |z-x| + |y|  \geq 1/C'_\eps.
\end{equation}
This proves (\ref{morse}).

In order to complete the  proof of (\ref{eq1.12.546}) we study in detail the characteristics of the flow generated by the vector field $\nu(x,y,z, \D_x, \D_y, \D_z)$. We first notice that the flow is complete since the gradient of $\Phi$ is Lipschitz.  Let $(x_0,y_0, z_0) \not\in \pi_{x,y,z}(\ccc + B(0,\eps))$. We denote by $(x(t), y(t), z(t))$ the integral curves of $\nu$ for  $t \in \R$ with
$x(0) = x_0$, $y(0) = y_0$, $z(0) = z_0$ and
$$
\left\{ \begin{array}{l}
\dot{x} = y \\
\dot{y} = -(\D_x V(x) - z) \\
\dot{z} = 0,
\end{array}
\right.
$$
and we have to study
\begin{equation} \label{f1}
 p_0 \sep{ \exp({t\nu}) (x_0, y_0, z_0)} = \frac{\gamma}{2} |x(t) - z(t)|^2 = \frac{\gamma}{2} |x(t) - z_0|^2.
\end{equation}
We split the study into two cases. Let $t\in [0, 1]$.

\paragraph{\it First case:} Suppose $|y_0| \geq |\D_xV(x_0) - z_0|$.
We write that
\begin{equation} \label{yy}
\begin{split}
y(t)  & = y_0 + \int_0^t \dot{y}(s) ds \\
& = y_0 - \int_0^t (\D_xV(x(s)) - z_0) ds \\
& = y_0 - t (\D_xV(x_0) - z_0) - \int_0^t \sep{ \int_0^s  V''(x(r))\dot{x}(r) dr} ds \\
& = y_0 - t (\D_xV(x_0) - z_0) - \int_0^t \sep{\int_0^s  V''(x(r))y(r) dr} ds.
\end{split}
\end{equation}
Since $V''$ is uniformly bounded, and denoting by $C_V \geq 1$ a corresponding bound we get
$$
\sup_{s \in [0,t]} |y(s) - y_0| \leq t |\D_xV(x_0) - z_0| + \frac{ t^2}{2} C_V \sup_{s \in [0,t]} |y(s)|.
$$
This implies that on $[0, t_V]$ with $t_V  \leq  1/C_V \leq 1/ C_V^{1/2}$ we have
$$
\frac{1}{2} \sup_{s \in [0,t]} |y(s) - y_0| \leq t |\D_xV(x_0) - z_0| + \frac{ t^2}{2} C_V |y_0|.
$$
Since by assumption $|y_0| \geq |\D_xV(x_0) - z_0|$ we get
$$
 \sup_{s \in [0,t]} |y(s) - y_0| \leq 2t |y_0| +  t^2 C_V |y_0| \leq 3 t |y_0|.
$$
We can then write that
\begin{equation}
\begin{split}
x(t) & = x_0 + \int_0^t \dot{x}(s) ds = x_0 + \int_0^t y(s) ds \\
    & = x_0 + t y_0 + \frac{3}{2}t^2 |y_0| \phi(t)
\end{split}
\end{equation}
with $|\phi(t)| \leq 1$ on $[0, t_V]$.
Recalling  (\ref{f1}) and using the triangular inequality, we get for $t\in [0, t_V]$ that
if $t_V \leq 1/4$
\begin{equation} \label{maxxx}
|x(t) - z_0|^2 \geq \max \set{  |x_0 - z_0| -2 t |y_0|,  \ \ \ \
                  \frac{t}{2} |y_0| - |x_0 - z_0|  }
\end{equation}
Let now $0 < \theta < t_V$. We split again the study into two parts:
\begin{enumerate}
\item If   $|x_0 - z_0| \geq \theta |y_0|$, we use the first expression in   (\ref{maxxx})
and we get that
$$
|x(t) - z_0|^2 \geq |x_0 - z_0|(1- 2t/ \theta) \geq \frac{1}{2} |x_0 - z_0| \ \ \textrm{ on } \ \ [0, \theta/4].
$$
\item If   $\theta |y_0| \geq |x_0 - z_0| $, we use the second expression in   (\ref{maxxx})
and we get that
$$
|x(t) - z_0|^2 \geq |y_0|(t/2-  \theta) \geq \frac{t_V}{8} |y_0| \ \ \textrm{ on } \ \ [t_V/2, t_V].
$$
if $\theta \leq t_V/8$.
\end{enumerate}
In all cases we get that there exists a constant   $c_V >0$
 depending only on $C_V$ such that
\begin{equation}  \label{estim1}
|x(t) - z_0|^2 \geq c_V \max \set{ |y_0|, |x_0 - z_0|, |\D_xV(x_0) - z_0|}
\end{equation}
on an interval of length at least $\theta$.

\paragraph{\it Second case:} Suppose $|y_0| \leq |\D_xV(x_0) - z_0|$.
As in (\ref{yy}) we can write that
\begin{equation}
\begin{split}
y(t)  = y_0 - t (\D_xV(x_0) - z_0) - \int_0^t \sep{\int_0^s  V''(x(r))y(r) dr} ds.
\end{split}
\end{equation}
Since $V''$ is uniformly bounded, and denoting again  by $C_V \geq 1$ a corresponding bound we get
$$
\sup_{s \in [0,t]} |y(s) - y_0 + s (\D_xV(x_0) - z_0)| \leq  \frac{ t^2}{2} C_V \sup_{s \in [0,t]} |y(s)|,
$$
so that
\begin{equation}
\begin{split}
 & \sup_{s \in [0,t]} |y(s) - y_0+ s (\D_xV(x_0) - z_0)| \\
 & \leq \frac{ t^2}{2} C_V  \sup_{s \in [0,t]} \abs{ y(s) - y_0+ s (\D_xV(x_0) - z_0)}
 + \frac{ t^2}{2} C_V \sep{ |y_0| + t |\D_xV(x_0) - z_0|}.
\end{split}
\end{equation}
This implies that on $[0, t_V]$ with again  $t_V  \leq  1/C_V \leq 1/ C_V^{1/2}$ we have
\begin{equation}
\begin{split}
\frac{1}{2} \sup_{s \in [0,t]} |y(s) - y_0+ s (\D_xV(x_0) - z_0)|
& \leq   \frac{ t^2}{2} C_V \sep{ |y_0| + t |\D_xV(x_0) - z_0|} \\
& \leq   C_V t^2   |\D_xV(x_0) - z_0|,
\end{split}
\end{equation}
since by assumption $|y_0| \leq |\D_xV(x_0) - z_0|$.  We therefore get
$$
 \sup_{s \in [0,t]} |y(s) - y_0 + s (\D_xV(x_0) - z_0)| \leq  2C_V t^2   |\D_xV(x_0) - z_0|.
$$
We can then write that
\begin{equation}
\begin{split}
x(t) & = x_0 + \int_0^t \dot{x}(s) ds = x_0 + \int_0^t y(s) ds \\
    & = x_0 + t y_0 - \frac{t^2}{2} (\D_xV(x_0) - z_0) + \frac{2}{3}  C_V t^3 |\D_xV(x_0) - z_0| \psi(t)
\end{split}
\end{equation}
with $|\psi(t)| \leq 1$ on $[0, t_V]$.
Recalling  (\ref{f1}) and using the triangular inequality, we get for $t\in [0, t_V]$ that
if $t_V \leq 3/ (8 C_V)$,
\begin{equation} \label{maxxx2}
\begin{split}
|x(t) - z_0|^2 \geq \max \left\{ \right. & |x_0 - z_0| - t |y_0| - t^2 |\D_xV(x_0) - z_0|, \\
  &  t |y_0|- |x_0 - z_0|  - t^2 |\D_xV(x_0) - z_0|, \\
   &                \frac{t^2}{4} |\D_xV(x_0) - z_0 |  - |x_0 - z_0| -  t|y_0| \left. \right\}
  \end{split}
\end{equation}
Let us now take again $0 < \theta < t_V$. We consider  three cases:
\begin{enumerate}
\item If   $|x_0 - z_0| \geq \max \set{ \theta |y_0|, 18 \theta^2 |\D_xV(x_0) - z_0 | }$, we use the first expression in   (\ref{maxxx2})
and we get that
$$
|x(t) - z_0|^2 \geq |x_0 - z_0|(1- t/ \theta -18 t^2/\theta^2) \geq \frac{1}{4} |x_0 - z_0| \ \ \textrm{ on } \ \ [0, \theta/6].
$$
\item If   $\theta |y_0| \geq \max \set{ |x_0 - z_0|,  18 \theta^2 |\D_xV(x_0) - z_0 |} $, we use the second expression in   (\ref{maxxx2})
and we get that
$$
|x(t) - z_0|^2 \geq |y_0|\sep{ t -  \theta - t^2/(18 \theta)} \geq \frac{\theta}{2} |y_0| \ \ \textrm{ on } \ \ [2 \theta, 3\theta].
$$
\item If   $18 \theta^2 |\D_xV(x_0) - z_0 | \geq \max \sep{ |x_0 - z_0|,  \theta |y_0|} $, we use the third expression in   (\ref{maxxx2})
and we get that
$$
|x(t) - z_0|^2 \geq |\D_xV(x_0) - z_0 |(\frac{1}{4} t^2 -  18 t \theta - 18 \theta^2) \geq \frac{t_V^2}{16} |\D_xV(x_0) - z_0 | \ \ \textrm{ on } \ \ [t_V/2, t_V].
$$
if $18 \theta^2 \leq t_V^2/16$.
\end{enumerate}
In all cases we get  that there exists a positive constant  $c_V'$ only depending on $C_V$ such that
\begin{equation}  \label{estim2}
 |x(t) - z_0| \geq c_V' \max \set{ |y_0|, |x_0 - z_0|, |\D_xV(x_0) - z_0|}
\end{equation}
on an interval of length at least $\theta$.

\paragraph{\it Conclusion} From (\ref{estim1}),  (\ref{estim2}) and using
(\ref{f1}), we get that in all cases
there exists an interval of length at least $\min (\theta, t_V/2)$ on which
\begeq
\begin{split}
p_0 \sep{ \exp({t\nu}) (x_0, y_0, z_0)} & =  \frac{\gamma}{2} |x(t) - z_0|^2  \\
& \geq \frac{\gamma}{2} \sep{ \min \sep{c_V', c_V} \max \set{ |y_0|, |x_0 - z_0|, |\D_xV(x_0) - z_0|}}^2.
\end{split}
\endeq
Since $(x_0,y_0, z_0) \not\in \pi_{xyz}(\ccc + B(0,\eps))$ and using (\ref{tuc3}) we get that the dynamical condition (\ref{eq1.12.546}) is fulfilled. The proof is complete. \fin

\end{document}